\documentclass[11pt,a4paper]{article}
\usepackage{amssymb}
\usepackage{amsmath}
\usepackage{amsthm}
\usepackage{verbatim}
\usepackage{pgf}
\usepackage{listings}
\usepackage{graphicx}
\usepackage{tabularx}
\usepackage{psfrag}
\usepackage{epstopdf}
\usepackage{multirow}
\usepackage{ifthen}
\usepackage{booktabs}
\usepackage{mathtools}
\usepackage{algorithm}
\usepackage{algorithmicx}
\usepackage{algpseudocode}
\usepackage{authblk}
\usepackage[justification=raggedleft,singlelinecheck=false]{caption}
\usepackage{geometry}
\numberwithin{equation}{section}

\newtheorem{thm}{Theorem}[section]  

\theoremstyle{definition}

\newtheorem{asn}[thm]{Assumption}
\newtheorem{lmm}[thm]{Lemma}
\newtheorem{prp}[thm]{Proposition}

\newtheorem{crl}[thm]{Corollary}

\newtheorem{rem}[thm]{Remark}

\theoremstyle{remark}
\newtheorem*{prf}{Proof}

\begin{document}

\title{Robust Inversion Methods for Aerosol Spectroscopy}

\author[a]{Graham Alldredge}
\author[b]{Tobias Kyrion}

\affil[a]{Center for Computational Engineering Science, RWTH Aachen University, Aachen, Germany}
\affil[b]{RWTH Aachen University, Aachen, Germany}

\renewcommand\Authands{ and }

\maketitle

\begin{abstract}
The Fast Aerosol Spectrometer (FASP) is a device for spectral aerosol measurements. Its purpose is to safely monitor the atmosphere inside a reactor containment. First we describe the FASP and explain its basic physical laws. Then we introduce our reconstruction methods for aerosol particle size distributions designed for the FASP. We extend known existence results for constrained Tikhonov regularization by uniqueness criteria and use those to generate reasonable models for the size distributions. We apply a Bayesian model-selection framework on these pre-generated models. We compare our algorithm with classical inversion methods using simulated measurements. We then extend our reconstruction algorithm for two-component aerosols, so that we can simultaneously retrieve their particle-size distributions and unknown volume fractions of their two components. Finally we present the results of a numerical study for the extended algorithm.
\end{abstract}

%

\section{The FASP measurement device}

The FASP is an optical measurement device for aerosol particle size distributions in rigid environments where the temperature may surpass $200^\circ\text{C}$ and the pressure 8 bar over atmospheric pressure, cf. \cite{kr15, kr14}. The aerosol particles themselves may be acidic as well. The FASP is split into a detector head and into a unit containing an evaluation computer and a light source with different adjustable light wavelengths. The sensitive evaluation and light source unit is connected with the robust detector head via two optical fibers.

The detector head is the only part of the FASP which extends into the containment with the aerosol to be measured and it consists of a pneumatically propelled tube. By moving the tube one can adjust a short or a long measurement path, where the two path lengths are $400$ and $800$ mm respectively. The sought-after aerosol particle size distributions are reconstructed from the light intensity loss on the gap distance between long and short path, so the FASP works in a similar way to a White cell. The detector head is equipped with two light detectors. The first one can receive light with wavelengths in the infrared domain from 0.8 - 3.4 $\mu$m, and the other one in the visible domain from 0.5 - 0.8 $\mu$m. 

\newpage

\begin{figure}[h!]
\centering
\includegraphics[width=1.0\textwidth]{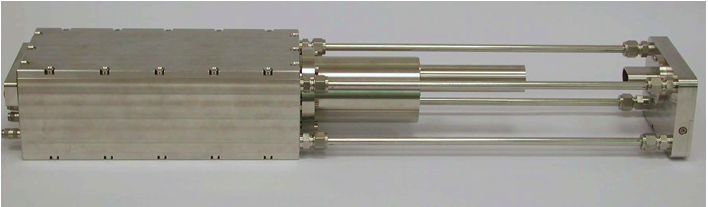}
\caption{The detector head with the movable tube (source: \cite{Sa11})}
\end{figure}

The ends of the optical fibers have to be floated with protective gas to shield them from harmful aerosol particles. These particle-free sections have to be subtracted from the actual geometric path lengths. This is not problematic since this does not change the gap distance. 

Let $l$ denote a current light wavelength used in a measurement, $G_{long}$ the geometric or unfloated long path and $G_{short}$ the geometric short path. The section floated with protective gas is labeled with $x$. Then the true path lengths are given by $L_{long} := G_{long} - x$ and $L_{short} := G_{short} - x.$

Let $M_{long}(l)$ and $M_{short}(l)$ be the measured intensities for long and short path, both perturbed by detector offsets $O_{long}(l)$ and $O_{short}(l)$ caused by ambient radiation. 

Then the intensities cleaned from the detector offsets are given by $I_{long}(l) := M_{long}(l) - O_{long}(l)$ and $I_{short}(l) := M_{short}(l) - O_{short}(l).$ 

According to the law of Beer-Lambert we have the relation
\begin{equation}
I_{long}(l) = I_{short}(l) \exp\left({\displaystyle{-\left(L_{long} - L_{short} \right) \int_{0}^{\infty} k(r,l) n(r) dr}}\right),
\end{equation}
where $n(r)$ is the sought-after unknown particle size distribution. The kernel function $k(r, l) := \pi r^{2}Q_{ext}(m_{med}(l), m_{part}(l), r, l)$ depends on both complex refractive indices $m_{med}(l)$ and $m_{part}(l)$ of the surrounding medium  and the scattering aerosol particles which depend on the wavelength $l$ of the incident light. The Mie extinction efficiency $Q_{ext}(m_{med}(l), m_{part}(l), r, l)$ of a spherical particle with radius $r$ illuminated by light with wavelength $l$ is derived from the general solution to the corresponding boundary value problem for Maxwell's equations and was first introduced in the pioneering article \cite{Mi08}. We adopt the numerical approximation of the Mie extinction efficiency in an absorbing medium from \cite{FS01}. From all of this follows
\begin{equation}
\int_{0}^{\infty} k(r,l) n(r) dr = e(l) \quad \text{with} \quad e(l) = - \frac{\log (I_{long}(l)) - \log(I_{short}(l))}{L_{long} - L_{short}}.
\label{operator_equation}
\end{equation}

\section{Modeling of FASP Measurement Data Inversions}
\label{motivation}

Let the measurement data $e(l)$ be an error-contaminated right-hand side for \eqref{operator_equation} and $(Kn)(l) := \int_{0}^{\infty} k(r,l) n(r) dr$ the compact linear operator with unbounded inverse which maps possible size distributions $n(r)$ to the left-hand side of \eqref{operator_equation}. We wish to reconstruct $n(r)$ from $e(l)$ by inverting the equation
\begin{equation}
Kn = e.
\label{operator_eq}
\end{equation}  
Here and in the following we omit the dependence on $r$ and $l$ for better readability. We assume that $e$ is given as a vector of finitely many independent Gaussian random variables $e_{i}$ with standard deviations $\sigma_{i}$ and means $\mu_{i}$, i.e.\ $e_{i} \sim \mathcal{N}(\mu_{i}, \sigma_{i}^{2})$. In the framework of Bayesian inference these are our \textit{observed random variables}. Now let $\boldsymbol{n} \in \mathbb{R}^{N}$ be a discrete approximation to $n$ and $\boldsymbol{K}_{N} \in \mathbb{R}^{N_{l} \times N}$ the \textit{kernel matrix} which correspondingly approximates the integral operator $K$. The details of these discretizations will be given in Section \ref{simulated_measurements}. We set up the covariance matrix $\boldsymbol{\Sigma_{\sigma}} = \mathrm{diag}(\sigma_{1}^{2}, ..., \sigma_{N_{l}}^{2})$. Then the \textit{observed model uncertainty} under the assumption $(\boldsymbol{K}_{N}\boldsymbol{n})_{i} = \mu_{i}$ obeys the probability distribution
\begin{equation}
p_{observed}(\boldsymbol{e} | \boldsymbol{n}) \propto \exp(-\textstyle{\frac{1}{2}}\|\boldsymbol{\Sigma_{\sigma}}^{-\frac{1}{2}}(\boldsymbol{K}_{N}\boldsymbol{n} - \boldsymbol{e})\|_{2}^{2}).
\end{equation}
After selecting a subjective \textit{prior distribution} $p_{prior}(\boldsymbol{n})$ which incorporates known \textit{a priori} information about $\boldsymbol{n}$ independent from the observed variable $\boldsymbol{e}$ we use Bayes' rule to obtain the \textit{posterior distribution} $p_{posterior}(\boldsymbol{n} | \boldsymbol{e})$ with
\begin{equation}
p_{posterior}(\boldsymbol{n} | \boldsymbol{e}) \propto p_{observed}(\boldsymbol{e} | \boldsymbol{n}) \times p_{prior}(\boldsymbol{n}).
\end{equation}
A more elaborate presentation of this Bayesian framework will be given in Section \ref{section_model_selection}. By applying a Tikhonov prior distribution
\begin{equation*}
p_{prior}(\boldsymbol{n}) \propto \exp(-\textstyle{\frac{1}{2}}\gamma \|\boldsymbol{n}\|_{2}^{2})I_{S}(\boldsymbol{n}),
\end{equation*}
where $\gamma \geq 0$ is a regularization parameter and $I_{S}(\boldsymbol{n})$ is the indicator function of the convex set
\begin{equation*}
S := \{\boldsymbol{n} \in \mathbb{R}^{N} | \; \boldsymbol{Cn} \leq \boldsymbol{b} \} \quad \text{with} \quad \boldsymbol{C} \in \mathbb{R}^{k \times N}, \; \boldsymbol{b} \in \mathbb{R}^{k},
\end{equation*}
we obtain the posterior distribution
\begin{equation}
p_{posterior}(\boldsymbol{n}) \propto \exp(-\textstyle{\frac{1}{2}}\|\boldsymbol{\Sigma_{\sigma}}^{-\frac{1}{2}}(\boldsymbol{K}_{N}\boldsymbol{n} - \boldsymbol{e})\|_{2}^{2} -\textstyle{\frac{1}{2}}\gamma \|\boldsymbol{n}\|_{2}^{2}))I_{S}(\boldsymbol{n}).
\label{post_prob}
\end{equation}

The quantity of interest $\boldsymbol{n}$ is estimated by computing the maximizer of the posterior distribution which is called the \textit{maximum a posteriori estimator} (MAP). It is obtained by solving the quadratic programming problem
\begin{equation}
\boldsymbol{n}_{MAP}^{\gamma} := \underset{\boldsymbol{n} \: \in \: \mathbb{R}^{N}}{\mathrm{argmin}} \; \textstyle{\frac{1}{2}}\|\boldsymbol{\Sigma_{\sigma}}^{-\frac{1}{2}}(\boldsymbol{K}_{N}\boldsymbol{n} - \boldsymbol{e})\|_{2}^{2} + \textstyle{\frac{1}{2}}\gamma \|\boldsymbol{n}\|_{2}^{2} \quad \text{s.t.} \quad \boldsymbol{Cn} \leq \boldsymbol{b}.
\label{MAP}
\end{equation}
Note that $\|\boldsymbol{\Sigma_{\sigma}}^{-\frac{1}{2}}(\boldsymbol{K}_{N}\boldsymbol{n} - \boldsymbol{e})\|_{2}^{2} \sim \chi^{2}(N_{l})$, which gives $\mathbb{E}(\|\boldsymbol{\Sigma_{\sigma}}^{-\frac{1}{2}}(\boldsymbol{K}_{N}\boldsymbol{n} - \boldsymbol{e})\|_{2}^{2}) = N_{l}$.

A classical residual-based inference method is the so-called \textit{discrepancy principle}.  After selecting a \textit{Morozov safety factor} $\tau$ the regularization parameter $\gamma$ is determined by demanding $\|\boldsymbol{\Sigma_{\sigma}}^{-\frac{1}{2}}(\boldsymbol{K}_{N}\boldsymbol{n} - \boldsymbol{e})\|_{2}^{2} = \tau N_{l}$. A common choice for the safety factor is $\tau = 1.1$. We will give a more thorough introduction to the discrepancy principle and some results on it in Section \ref{constrained_Tikhonov}. 

\textit{Monte Carlo methods} offer another way to evaluate the posterior distribution, cf. \cite{Chi96}. The advantage of Monte Carlo methods is that they take more of the statistical behavior of the observed measurement noise into account because all possible solutions with nonnegliglible posterior probabibility are sampled and contribute to the inference result. However these methods require a lot of computational resources, which we cannot afford because our application requires that one FASP measurement data inversion must be completed in under thirty seconds using a regular notebook. 

The discrepancy principle gets along with much less computational effort, but it does not take into account the specific shape of the distribution of the observed measurement noise. It does not explore the posterior distribution thoroughly and might give unreasonable results because of this.

In our hybrid approach we combine the advantages of both methods. We review Tikhonov regularization under linear constraints and derive conditions for the existence of a bijection between the regularization parameter and the residual. If these conditions are fulfilled, we can propose a set of regularization parameters obtained with the discrepancy principle using a set of Morozov safety factors corresponding to high-probability values of the weighted norm of the residual, $\|\boldsymbol{\Sigma_{\sigma}}^{-\frac{1}{2}}(\boldsymbol{K}_{N}\boldsymbol{n} - \boldsymbol{e})\|_{2}^{2} \sim \chi^{2}(N_{l})$. After this a Bayesian model-comparison procedure is applied to these reconstructions, and we rank them according to their posterior probabilities.

We show with numerical simulations that our method satisfies the demands on runtime and accuracy and that it is superior to existing inversion methods based on classical model-selection approaches. In the last section we extend our method to investigate two-component aerosols.

\section{Tikhonov Regularization under Linear Constraints}
\label{constrained_Tikhonov}

Computing the maximum a posteriori estimator leads to a quadratic programming problem of the form
\begin{equation}
\boldsymbol{n}_{\gamma} := \underset{\boldsymbol{n} \: \in \: \mathbb{R}^{N}}{\mathrm{argmin}} \; \textstyle{\frac{1}{2}}\|\boldsymbol{Kn} - \boldsymbol{r}\|_{2}^{2} + \textstyle{\frac{1}{2}}\gamma\|\boldsymbol{n}\|_{2}^{2} \quad \text{s.t.} \quad \boldsymbol{Cn} \leq \boldsymbol{b},
\label{quadprog}
\end{equation}
with $\boldsymbol{K} := \boldsymbol{\Sigma_{\sigma}}^{-\frac{1}{2}}\boldsymbol{K}_{N}$ and $ \boldsymbol{r} := \boldsymbol{\Sigma_{\sigma}}^{-\frac{1}{2}}\boldsymbol{e}$. The function to be minimized is known as the Tikhonov functional. 

It is proved in \cite{EHN96} that the residual of the Tikhonov-regularized solution under linear constraints decreases monotonically with the regularization parameter $\gamma$. To the best of our knowledge conditions for \textit{strict} monotonicity have not been found yet, so we derive some in the following. The advantage of having a strictly monotonic relation between regularization parameter and residual is that it gives a bijection. Thus we can then identify any regularization parameter $\gamma$ from the range $[0, \infty)$ with a unique residual value $\|\boldsymbol{Kn}_{\gamma} - \boldsymbol{r}\|_{2}^{2}$ from the range $[\|\boldsymbol{Kn}_{0} - \boldsymbol{r}\|_{2}^{2}, \|\boldsymbol{Kn}_{\infty} - \boldsymbol{r}\|_{2}^{2})$. Here 
\begin{equation}
\label{minimum_norm}
\boldsymbol{n}_{\infty} := \underset{\boldsymbol{n} \: \in \: \mathbb{R}^{N}}{\mathrm{argmin}} \; \textstyle{\frac{1}{2}}\|\boldsymbol{n}\|_{2}^{2} \quad \text{s.t.} \quad \boldsymbol{Cn} \leq \boldsymbol{b}
\end{equation}
is the \textit{minimum norm element}. As shown in \cite{Neub88a} there holds $\lim_{\gamma \to \infty} \boldsymbol{n}_{\gamma} = \boldsymbol{n}_{\infty}$. When our monotonicity conditions are satisfied, we obtain a set of distinct regularization parameters by proposing a set of distinct residual values from the range $[\|\boldsymbol{Kn}_{0} - \boldsymbol{r}\|_{2}^{2}, \|\boldsymbol{Kn}_{\infty} - \boldsymbol{r}\|_{2}^{2})$. The disadvantageous case of multiple prior distributions corresponding to the same residual value can therefore not occur. Note that in practice the cases $\gamma = 0$ and $\gamma = \infty$ are inadmissible, since then the Tikhonov prior distribution is improper or degenerates to a point mass, so we always restrict ourselves to a finite range $(0, \gamma_{max}]$ with $\gamma_{max} < \infty $.

\subsection{Necessary Conditions for Strict Monotonicity}

The following theorem shows that $\boldsymbol{n}_{\alpha} \neq \boldsymbol{n}_{\beta}$ for all $\alpha > \beta$ is the only necessary condition needed for strict monotonicity.

\begin{lmm}
\label{monotonicity}
\textit{Let $\alpha > \beta \geq 0$ be arbitrary and $\boldsymbol{n}_{\alpha}$ and $\boldsymbol{n}_{\beta}$ the solutions of \eqref{quadprog} for $\gamma = \alpha$ and $\gamma = \beta$ respectively. If there holds $\boldsymbol{n}_{\alpha} \neq \boldsymbol{n}_{\beta}$ for all $\alpha > \beta$, then the residual $\|\boldsymbol{Kn}_{\gamma} - \boldsymbol{r}\|_2$ is strictly increasing for growing $\gamma$.}
\end{lmm}

\begin{prf}
From the first-order necessary Karush-Kuhn-Tucker conditions for the problem \eqref{quadprog} we have that for each $\gamma$ there exists a vector $\boldsymbol{q}_{\gamma} \in \mathbb{R}^{k}$ with
\begin{align}
\boldsymbol{K}^{T}\boldsymbol{K}\boldsymbol{n}_{\gamma} - \boldsymbol{K}^{T}\boldsymbol{r} + \gamma\boldsymbol{n}_{\gamma} + \boldsymbol{C}^{T}\boldsymbol{q}_{\gamma} & = 0 \label{KKT_lin_constr_1}\\
\boldsymbol{Cn}_{\gamma} & \leq \boldsymbol{b} \label{KKT_lin_constr_3}\\
\mathrm{diag}\left(\boldsymbol{q}_{\gamma}\right)\left(\boldsymbol{Cn}_{\gamma} - \boldsymbol{b}\right) & = 0 \label{KKT_lin_constr_4}\\
\boldsymbol{q}_{\gamma} & \geq 0. \label{KKT_lin_constr_5} 
\end{align}

We define the difference vector
\begin{equation*}
\boldsymbol{x} := \boldsymbol{n}_{\beta} - \boldsymbol{n}_{\alpha}
\end{equation*}
and subtract \eqref{KKT_lin_constr_1} for $\gamma = \alpha$ with the same equation for $\gamma = \beta$ to get
\begin{equation}
\boldsymbol{K}^{T}\boldsymbol{Kx} + \beta\boldsymbol{n}_{\beta} - \alpha\boldsymbol{n}_{\alpha} +  \boldsymbol{C}^{T}\left(\boldsymbol{q}_{\beta} - \boldsymbol{q}_{\alpha}\right) = 0.
\label{KKT-diff}
\end{equation}
Taking the scalar product of \eqref{KKT-diff} with $\boldsymbol{n}_{\alpha}$ and then with $\boldsymbol{n}_{\beta}$ gives
\begin{align*}
\left\langle \boldsymbol{n}_{\alpha}, \boldsymbol{K}^{T}\boldsymbol{Kx} + \beta\boldsymbol{n}_{\beta} - \alpha\boldsymbol{n}_{\alpha} + \boldsymbol{C}^{T}\left(\boldsymbol{q}_{\beta} - \boldsymbol{q}_{\alpha}\right) \right\rangle & = 0\\ 
\text{and} \quad \left\langle \boldsymbol{n}_{\beta}, \boldsymbol{K}^{T}\boldsymbol{Kx} + \beta\boldsymbol{n}_{\beta} - \alpha\boldsymbol{n}_{\alpha} + \boldsymbol{C}^{T}\left(\boldsymbol{q}_{\beta} - \boldsymbol{q}_{\alpha}\right) \right\rangle & = 0.
\end{align*}
Our next step is to add $(\alpha - \beta)\left\langle \boldsymbol{n}_{\alpha}, \boldsymbol{n}_{\alpha} \right\rangle$
on both sides of the first relation and analogously \newline $(\alpha - \beta)\left\langle \boldsymbol{n}_{\beta}, \boldsymbol{n}_{\beta} \right\rangle$ on both sides of the latter relation, which results in
\begin{align*}
& \left\langle \boldsymbol{n}_{\alpha}, \left(\boldsymbol{K}^{T}\boldsymbol{K} + \beta \boldsymbol{I}\right)\boldsymbol{x} +  \boldsymbol{C}^{T}\left(\boldsymbol{q}_{\beta} - \boldsymbol{q}_{\alpha}\right) \right\rangle = (\alpha - \beta)\left\langle \boldsymbol{n}_{\alpha}, \boldsymbol{n}_{\alpha} \right\rangle\\
\text{and} \qquad & \left\langle \boldsymbol{n}_{\beta}, \left(\boldsymbol{K}^{T}\boldsymbol{K} + \alpha \boldsymbol{I}\right)\boldsymbol{x} + \boldsymbol{C}^{T}\left(\boldsymbol{q}_{\beta} - \boldsymbol{q}_{\alpha}\right) \right\rangle = (\alpha - \beta)\left\langle \boldsymbol{n}_{\beta}, \boldsymbol{n}_{\beta} \right\rangle . 
\end{align*}
Taking the difference of these two equations gives
\begin{align*}
& \;\: (\alpha - \beta)\big(\left\langle \boldsymbol{n}_{\beta}, \boldsymbol{n}_{\beta} \right\rangle - \left\langle \boldsymbol{n}_{\alpha}, \boldsymbol{n}_{\alpha} \right\rangle\big) \\
= & \; \left\langle \boldsymbol{x}, \left(\boldsymbol{K}^{T}\boldsymbol{K} + \alpha \boldsymbol{I}\right)\boldsymbol{n}_{\beta} \right\rangle - \left\langle \boldsymbol{x}, \left(\boldsymbol{K}^{T}\boldsymbol{K} + \beta \boldsymbol{I}\right)\boldsymbol{n}_{\alpha}\right\rangle + \left\langle \boldsymbol{x}, \boldsymbol{C}^{T}\left(\boldsymbol{q}_{\beta} - \boldsymbol{q}_{\alpha}\right) \right\rangle .
\end{align*}

On the one hand this implies
\begin{align*}
&  \;\: (\alpha - \beta)\big(\left\langle \boldsymbol{n}_{\beta}, \boldsymbol{n}_{\beta} \right\rangle - \left\langle \boldsymbol{n}_{\alpha}, \boldsymbol{n}_{\alpha} \right\rangle\big) \\
= & \; \left\langle \boldsymbol{x}, \left(\boldsymbol{K}^{T}\boldsymbol{K} + \beta \boldsymbol{I}\right)\boldsymbol{n}_{\beta} \right\rangle + (\alpha - \beta)\left\langle \boldsymbol{x}, \boldsymbol{n}_{\beta}\right\rangle - \left\langle \boldsymbol{x}, \left(\boldsymbol{K}^{T}\boldsymbol{K} + \beta \boldsymbol{I}\right)\boldsymbol{n}_{\alpha}\right\rangle + \left\langle \boldsymbol{x}, \boldsymbol{C}^{T}\left(\boldsymbol{q}_{\beta} - \boldsymbol{q}_{\alpha}\right) \right\rangle \\
= & \; \left\langle \boldsymbol{x}, \left(\boldsymbol{K}^{T}\boldsymbol{K} + \beta \boldsymbol{I}\right)\boldsymbol{x} \right\rangle + (\alpha - \beta)\left\langle \boldsymbol{x}, \boldsymbol{n}_{\beta}\right\rangle + \left\langle \boldsymbol{x}, \boldsymbol{C}^{T}\left(\boldsymbol{q}_{\beta} - \boldsymbol{q}_{\alpha}\right) \right\rangle,
\end{align*}
while on the other hand
\begin{align*}
& \;\: (\alpha - \beta)\big(\left\langle \boldsymbol{n}_{\beta}, \boldsymbol{n}_{\beta} \right\rangle - \left\langle \boldsymbol{n}_{\alpha}, \boldsymbol{n}_{\alpha} \right\rangle\big) \\
= & \; \left\langle \boldsymbol{x}, \left(\boldsymbol{K}^{T}\boldsymbol{K} + \alpha \boldsymbol{I}\right)\boldsymbol{n}_{\beta} \right\rangle  - \left\langle \boldsymbol{x}, \left(\boldsymbol{K}^{T}\boldsymbol{K} + \alpha \boldsymbol{I}\right)\boldsymbol{n}_{\alpha}\right\rangle - (\beta - \alpha)\left\langle \boldsymbol{x}, \boldsymbol{n}_{\alpha}\right\rangle + \left\langle \boldsymbol{x}, \boldsymbol{C}^{T}\left(\boldsymbol{q}_{\beta} - \boldsymbol{q}_{\alpha}\right) \right\rangle \\
= & \; \left\langle \boldsymbol{x}, \left(\boldsymbol{K}^{T}\boldsymbol{K} + \alpha \boldsymbol{I}\right)\boldsymbol{x} \right\rangle + (\alpha - \beta)\left\langle \boldsymbol{x}, \boldsymbol{n}_{\alpha}\right\rangle + \left\langle \boldsymbol{x}, \boldsymbol{C}^{T}\left(\boldsymbol{q}_{\beta} - \boldsymbol{q}_{\alpha}\right) \right\rangle
\end{align*}
holds.
Adding these gives
\begin{align*}
& \;\: 2(\alpha - \beta)\big(\left\langle \boldsymbol{n}_{\beta}, \boldsymbol{n}_{\beta} \right\rangle - \left\langle \boldsymbol{n}_{\alpha}, \boldsymbol{n}_{\alpha} \right\rangle\big) \\
= & \; \left\langle \boldsymbol{x}, \left(2\boldsymbol{K}^{T}\boldsymbol{K} + (\alpha + \beta) \boldsymbol{I}\right)\boldsymbol{x} \right\rangle + (\alpha - \beta)\big(\left\langle \boldsymbol{n}_{\beta}, \boldsymbol{n}_{\beta} \right\rangle - \left\langle \boldsymbol{n}_{\alpha}, \boldsymbol{n}_{\alpha} \right\rangle\big) + 2\left\langle \boldsymbol{x}, \boldsymbol{C}^{T}\left(\boldsymbol{q}_{\beta} - \boldsymbol{q}_{\alpha}\right) \right\rangle,
\end{align*}
and finally we arrive at
\begin{equation}
\begin{split}
& \;\: (\alpha - \beta)\big(\left\langle \boldsymbol{n}_{\beta}, \boldsymbol{n}_{\beta} \right\rangle - \left\langle \boldsymbol{n}_{\alpha}, \boldsymbol{n}_{\alpha} \right\rangle\big) \\
= & \; \left\langle \boldsymbol{n}_{\beta} - \boldsymbol{n}_{\alpha}, \left(2\boldsymbol{K}^{T}\boldsymbol{K} + (\alpha + \beta) \boldsymbol{I}\right)\left(\boldsymbol{n}_{\beta} - \boldsymbol{n}_{\alpha}\right)\right\rangle + 2\left\langle \boldsymbol{C}\left(\boldsymbol{n}_{\beta} - \boldsymbol{n}_{\alpha}\right), \boldsymbol{q}_{\beta} - \boldsymbol{q}_{\alpha} \right\rangle.
\end{split}
\label{main_eq}
\end{equation}

Now we consider the term $\left\langle \boldsymbol{C}\left(\boldsymbol{n}_{\beta} - \boldsymbol{n}_{\alpha}\right), \boldsymbol{q}_{\beta} - \boldsymbol{q}_{\alpha} \right\rangle$. 
The following four cases can occur:

\begin{table}[h!]
\centering
\begin{tabular}{||m{1.2cm}|m{1.2cm}||m{1.1cm}|m{1.1cm}|m{1.5cm}||m{1.1cm}|m{1.1cm}|m{1.1cm}||}
\hline
\multicolumn{2}{||c||}{$i$-th constraint} & \centering \multirow{4}{*}{$\left(\boldsymbol{C}\boldsymbol{n}_{\beta}\right)_{i}$\vspace{-0.2cm}} & \centering \multirow{4}{*}{$\left(\boldsymbol{C}\boldsymbol{n}_{\alpha}\right)_{i}$\vspace{-0.2cm}} & \centering \multirow{4}{*}{\hspace{0.36cm}$\left(\boldsymbol{C}\boldsymbol{n}_{\beta}\right)_{i}$ \vspace{0.3cm}} & \centering \multirow{4}{*}{$\left(\boldsymbol{q}_{\beta}\right)_{i}$\vspace{-0.2cm}} & \centering \multirow{4}{*}{$\left(\boldsymbol{q}_{\alpha}\right)_{i}$\vspace{-0.2cm}} & \centering \multirow{4}{*}{\hspace{0.36cm}$\left(\boldsymbol{q}_{\beta}\right)_{i}$ \vspace{0.3cm}}\tabularnewline
\multicolumn{2}{||c||}{for the Tikhonov} & & & & & & \tabularnewline
\multicolumn{2}{||c||}{functional for} & & &  {$- \left(\boldsymbol{C}\boldsymbol{n}_{\alpha}\right)_{i}$ \vspace{-0.8cm}} & & & {$- \left(\boldsymbol{q}_{\alpha}\right)_{i}$ \vspace{-0.8cm}} \tabularnewline
\cline{1-2}
\centering $\gamma = \beta$ & \centering $\gamma = \alpha$ & & & & & & \tabularnewline
\hline
\hline
\centering active & \centering active & \centering $= \left(\boldsymbol{b}\right)_{i}$ & \centering $= \left(\boldsymbol{b}\right)_{i}$ & \centering $= 0$ & \centering $\geq 0$ & \centering $\geq 0$ & \centering void\tabularnewline
\hline
\centering inactive &  \centering active & \centering $< \left(\boldsymbol{b}\right)_{i}$ & \centering $= \left(\boldsymbol{b}\right)_{i}$ & \centering $< 0$ & \centering $= 0$ & \centering $\geq 0$ & \centering $\leq 0$ \tabularnewline
\hline
\centering active &  \centering inactive & \centering $= \left(\boldsymbol{b}\right)_{i}$  & \centering $< \left(\boldsymbol{b}\right)_{i}$ & \centering $> 0$ & \centering $\geq 0$ & \centering $= 0$ & \centering $\geq 0$ \tabularnewline
\hline
\centering inactive &  \centering inactive & \centering $< \left(\boldsymbol{b}\right)_{i}$ & \centering $< \left(\boldsymbol{b}\right)_{i}$ & \centering void & \centering $= 0$ & \centering $= 0$ & \centering $= 0$ \tabularnewline
\hline
\end{tabular}
\end{table}

From this we see that all components of the vector 
\begin{equation*}
\mathrm{diag}\big( \boldsymbol{C}\left(\boldsymbol{n}_{\beta} - \boldsymbol{n}_{\alpha}\right)\big)\left(\boldsymbol{q}_{\beta} - \boldsymbol{q}_{\alpha}\right)
\end{equation*}
are nonnegative, and so $\left\langle \boldsymbol{C}\left(\boldsymbol{n}_{\beta} - \boldsymbol{n}_{\alpha}\right), \boldsymbol{q}_{\beta} - \boldsymbol{q}_{\alpha} \right\rangle \geq 0$. Under the assumption $\boldsymbol{n}_{\alpha} \neq \boldsymbol{n}_{\beta}$ we have $\boldsymbol{x} \neq 0$, and since the matrix $2\boldsymbol{K}^{T}\boldsymbol{K} + (\alpha + \beta) \boldsymbol{I}$ is positive definite we finally conclude with \eqref{main_eq} that
\begin{equation*}
(\alpha - \beta)\big(\left\langle \boldsymbol{n}_{\beta}, \boldsymbol{n}_{\beta} \right\rangle - \left\langle \boldsymbol{n}_{\alpha}, \boldsymbol{n}_{\alpha} \right\rangle\big) > 0, 
\end{equation*}
which is equivalent to $\|\boldsymbol{n}_{\beta}\|_{2}^{2} > \|\boldsymbol{n}_{\alpha}\|_{2}^{2}$.

We proceed then with
\begin{equation}
\|\boldsymbol{Kn}_{\alpha} - \boldsymbol{r}\|_{2}^{2} - \|\boldsymbol{Kn}_{\beta} - \boldsymbol{r}\|_{2}^{2} = \left\langle \boldsymbol{x}, \boldsymbol{K}^{T}\boldsymbol{Kx} \right\rangle + 2\left\langle -\boldsymbol{x}, \boldsymbol{K}^{T}\boldsymbol{Kn}_{\beta} - \boldsymbol{K}^{T}\boldsymbol{r} + \beta\boldsymbol{n}_{\beta} \right\rangle + 2\beta\left\langle \boldsymbol{x}, \boldsymbol{n}_{\beta} \right\rangle 
\label{res_eq}
\end{equation}  
by using $\boldsymbol{n}_{\alpha} = \boldsymbol{n}_{\beta} - \boldsymbol{x}$. The variational inequality for the Tikhonov functional for $\gamma = \beta$ yields $\left\langle -\boldsymbol{x}, \boldsymbol{K}^{T}\boldsymbol{Kn}_{\beta} - \boldsymbol{K}^{T}\boldsymbol{r} + \beta\boldsymbol{n}_{\beta} \right\rangle \geq 0$. Moreover we have
\begin{align*}
\left\langle \boldsymbol{x}, \boldsymbol{n}_{\beta} \right\rangle \; = \; \left\langle \boldsymbol{n}_{\beta} - \boldsymbol{n}_{\alpha}, \boldsymbol{n}_{\beta} \right\rangle
\; \geq \; \|\boldsymbol{n}_{\beta}\|_{2}^{2} - \|\boldsymbol{n}_{\alpha}\|_{2}\|\boldsymbol{n}_{\beta}\|_{2}
\; > \; 0.
\end{align*}
In summary we have shown $\|\boldsymbol{Kn}_{\alpha} - \boldsymbol{r}\|_{2}^{2} > \|\boldsymbol{Kn}_{\beta} - \boldsymbol{r}\|_{2}^{2}.$ \hfill $\square$ 
\end{prf}

\begin{rem}
\label{res_stuck}
From \eqref{res_eq} follows that all $\boldsymbol{n}_{\gamma}$ with $\|\boldsymbol{Kn}_{\gamma} - \boldsymbol{r}\|_2 = \tau$ for an arbitrary but fixed $\tau$ must coincide. This means in other words that if the residual of the regularized solutions `gets stuck' at some value $\tau$, the solutions $\boldsymbol{n}_{\gamma}$ are constant for these values of $\gamma$. In the next section we derive conditions which prevent this case.  
\end{rem}

\subsection{Sufficient Conditions for Strict Monotonicity}

In this section we derive sufficient conditions for $\boldsymbol{n}_{\alpha} \neq \boldsymbol{n}_{\beta}$ for $\alpha > \beta$, hence by Lemma \ref{monotonicity} for strict monotonicity. In particular we focus on constraints of the form $\boldsymbol{Cn} \geq 0$ with $\boldsymbol{C} \in \mathbb{R}^{k \times N}$ and $k \leq N$, i.e.\ on generalized nonnegativity constraints. For this specific type of constraints we have that for the minimum norm solution $\boldsymbol{n}_{\infty}$ defined in \eqref{minimum_norm} that $\boldsymbol{n}_{\infty} \equiv 0$ holds, which gives according to \cite{Neub88a} the relation $\|\boldsymbol{Kn}_{\alpha} - \boldsymbol{r}\|_{2} \leq \|\boldsymbol{r}\|_{2}$ for all $\alpha \geq 0$. We carry out our following considerations under following important assumption.
\begin{asn}
\label{noise_assumption} 
\textit{The regularization parameter $\alpha$ is selected in such a way that for the regularized residual the inequality $\|\boldsymbol{Kn}_{\alpha} - \boldsymbol{r}\|_{2} \leq c\delta$ holds, where $c > 0$ is a fixed constant and $\delta > 0$ the noise level so that the relation $c\delta < \|\boldsymbol{r}\|_{2}$ is satisfied.}
\end{asn}
Under this assumption we have the strict inequality $\|\boldsymbol{Kn}_{\alpha} - \boldsymbol{r}\|_{2} < \|\boldsymbol{r}\|_{2}$ for all $\alpha \in [0, \infty)$.
 
\begin{rem}
If we assume $\boldsymbol{Kn}_{\alpha} = \boldsymbol{r}_{true}$, where $\boldsymbol{r}_{true}$ is the ``true'' data vector, we can rewrite the first part of above assumption as $\|\boldsymbol{r} - \boldsymbol{r}_{true}\|_{2} < c\delta$. This is a standard assumption made in inverse problems literature, cf. \cite{Neub88a}. 
\end{rem}

\begin{thm}
\label{monotonicity_cond}
\textit{Let $\boldsymbol{n}_{\alpha}$ be given by}
\begin{equation}
 \boldsymbol{n}_{\alpha} := \underset{\boldsymbol{n} \: \in \: \mathbb{R}^{N}}{\mathrm{argmin}} \; \textstyle{\frac{1}{2}}\|\boldsymbol{Kn} - \boldsymbol{r}\|_{2}^{2} + \textstyle{\frac{1}{2}}\alpha\|\boldsymbol{n}\|_{2}^{2} \quad \text{s.t.} \quad -\boldsymbol{Cn} \leq 0,
\label{quadprog_nng}
\end{equation}
with $\boldsymbol{C} \in \mathbb{R}^{k \times N}$ having full row rank $k \leq N$. If $\|\boldsymbol{Kn}_{\alpha} - \boldsymbol{r}\|_{2} < \|\boldsymbol{r}\|_{2}$, or equivalently $\boldsymbol{n}_{\alpha} \neq 0$ for all $\alpha \in [0, \infty)$ according to Lemma \ref{monotonicity} and Remark \ref{res_stuck}, then we have $\boldsymbol{n}_{\alpha} \neq \boldsymbol{n}_{\beta}$ for all $\alpha > \beta$.
\end{thm}

\begin{prf}
Let $\alpha > \beta$. Let $\boldsymbol{C}_{act}^{\alpha}$ denote the submatrix of $\boldsymbol{C}$ with active constraints in \eqref{quadprog_nng} for the regularization parameter $\alpha$.

We first consider the case $\boldsymbol{C}_{act}^{\alpha} \neq \boldsymbol{C}_{act}^{\beta}$. We obtain $\boldsymbol{C}(\boldsymbol{n}_{\alpha} - \boldsymbol{n}_{\beta}) \neq 0$, i.e.\ $\boldsymbol{n}_{\alpha} - \boldsymbol{n}_{\beta} \notin \mathrm{ker}(\boldsymbol{C})$. This gives directly $\boldsymbol{n}_{\alpha} - \boldsymbol{n}_{\beta} \neq 0$. 

Now we turn to the case $\boldsymbol{C}_{act}^{\alpha} = \boldsymbol{C}_{act}^{\beta}$. The first-order necessary conditions for a minimizer in \eqref{quadprog_nng} are given by
\begin{equation}
\boldsymbol{K}^{T}\boldsymbol{Kn} - \boldsymbol{K}^{T}\boldsymbol{r} + \alpha\boldsymbol{n} - \boldsymbol{C}^{T}\boldsymbol{q}_{\alpha} = 0,
\label{KKT_nng}
\end{equation}
where $\boldsymbol{q}_{\alpha} \geq 0$. Let us assume $\boldsymbol{n}_{\alpha} = \boldsymbol{n}_{\beta}$. Then taking the difference of \eqref{KKT_nng} for the parameters $\alpha$ and $\beta$ yields
\begin{equation}
(\alpha - \beta)\boldsymbol{n}_{\alpha} - \boldsymbol{C}^{T}(\boldsymbol{q}_{\alpha} - \boldsymbol{q}_{\beta}) = 0.
\label{diff_eq}
\end{equation}
Let us first consider the subcase that none of the constraints is active. Then we have $\boldsymbol{q}_{\alpha} = \boldsymbol{q}_{\beta} = 0$, which implies $(\alpha - \beta)\boldsymbol{n}_{\alpha} = 0$. This contradicts $\boldsymbol{n}_{\alpha} \neq 0$, so we must have $\boldsymbol{n}_{\alpha} \neq \boldsymbol{n}_{\beta}$. Now we turn to the subcase that at least one constraint is active. Let $\boldsymbol{q}_{act}^{\alpha}$ and $\boldsymbol{q}_{act}^{\beta}$ be the subvectors of $\boldsymbol{q}_{\alpha}$ and $\boldsymbol{q}_{\beta}$ corresponding to active constraints. Then we can rewrite the last equation as
\begin{equation*}
(\alpha - \beta)\boldsymbol{n}_{\alpha} - {\boldsymbol{C}_{act}^{\alpha}}^{T}(\boldsymbol{q}_{act}^{\alpha} - \boldsymbol{q}_{act}^{\beta}) = 0,
\end{equation*}
where we remember that $\boldsymbol{C}_{act}^{\alpha}$ is obtained from $\boldsymbol{C}$ by canceling its $i$-th row when the constraint $-(\boldsymbol{C}\boldsymbol{n})_{i} \leq 0$ is inactive and thus $(\boldsymbol{q}_{\alpha})_{i} = (\boldsymbol{q}_{\beta})_{i} = 0$ holds. Our next step is to multiply this equation from the left with $\boldsymbol{C}_{act}^{\alpha}$. By construction of $\boldsymbol{C}_{act}^{\alpha}$ we have $\boldsymbol{C}_{act}^{\alpha}\boldsymbol{n}_{\alpha} = 0$ and therefore
\begin{equation*}
- \; \boldsymbol{C}_{act}^{\alpha}{\boldsymbol{C}_{act}^{\alpha}}^{T}(\boldsymbol{q}_{act}^{\alpha} - \boldsymbol{q}_{act}^{\beta}) = 0.
\end{equation*}
Since $\boldsymbol{C}_{act}^{\alpha}{\boldsymbol{C}_{act}^{\alpha}}^{T}$ has full rank, this implies $\boldsymbol{q}_{act}^{\alpha} = \boldsymbol{q}_{act}^{\beta}$ and hence $\boldsymbol{q}_{\alpha} = \boldsymbol{q}_{\beta}$. Inserting this finding back into \eqref{diff_eq} gives $\boldsymbol{n}_{\alpha} = 0$, which contradicts our assumption $\boldsymbol{n}_{\alpha} \neq 0$ for all $\alpha \in [0, \infty)$. Therefore we must also have $\boldsymbol{n}_{\alpha} \neq \boldsymbol{n}_{\beta}$ in this subcase.
\hfill $\square$
\end{prf}

\subsection{The Discrepancy Principle}
\label{discrepancy_principle}

With the next Theorem we summarize our previous results.
\begin{thm}
\label{discrepancy_principle}
Let the conditions of Theorem \ref{monotonicity_cond} be fulfilled and let $\boldsymbol{n}_{\infty}$ be the minimum norm solution defined in \eqref{minimum_norm}. Define $r_{0} := \|\boldsymbol{Kn}_{0} - \boldsymbol{r}\|_2$ and $r_{\infty} :=  \|\boldsymbol{Kn}_{\infty} - \boldsymbol{r}\|_2$. Then there exist for any $\tau$ from $[r_{0}, r_{\infty})$ a unique $\gamma$ from $[0, \infty)$ such that $\|\boldsymbol{Kn}_{\gamma} - \boldsymbol{r}\|_2 = \tau$. The residual grows strictly monotonically with $\gamma$. 
\end{thm}
\hfill $\square$

\begin{rem}
The discrepancy principle carries directly over to generalized Tikhonov regularization, where the prior distribution is given by
\begin{equation*}
p_{prior}(\boldsymbol{n}) \propto \exp(-\textstyle{\frac{1}{2}}\gamma \boldsymbol{n}^{T}\boldsymbol{R}\boldsymbol{n})I_{S}(\boldsymbol{n}),
\end{equation*}
where $\boldsymbol{R}$ is a positive definite regularization matrix and $I_{S}(\boldsymbol{n})$ is the indicator function of $S = \{\boldsymbol{n} \in \mathbb{R}^{N} | \; - \boldsymbol{Cn} \leq 0 \}$. Here we have to solve the quadratic programming problem
\begin{equation*}
\min_{\boldsymbol{n} \in \mathbb{R}^{N}} \textstyle{\frac{1}{2}}\|\boldsymbol{Kn} - \boldsymbol{r}\|_{2}^{2} + \textstyle{\frac{1}{2}}\gamma \boldsymbol{n}^{T}\boldsymbol{R}\boldsymbol{n} \quad \text{s.t.} \quad - \boldsymbol{Cn} \leq 0.
\end{equation*}
Let $\boldsymbol{R} = \boldsymbol{U}^{T}\boldsymbol{U}$ the Cholesky decomposition. Then the substitution $\boldsymbol{n} = \boldsymbol{U}^{-1}\boldsymbol{v}$ transforms the above quadratic programming problem into the standard form \eqref{quadprog}.
\end{rem}

\subsection{Convergence Analysis}
\label{convergence1}

At this point we review some classical convergence criteria for parameter-choice strategies for Tikhonov regularization under linear constraints. With convergence we mean that the regularized reconstructions approach the true solution of the noise-free linear inverse problem as the noise level goes to $0$. We decompose the noisy data vector $\boldsymbol{r}$ into
\begin{align*}
\boldsymbol{r} & = \boldsymbol{\Sigma_{\sigma}}^{-\frac{1}{2}}(\boldsymbol{e}_{true} + \boldsymbol{\delta}) \\
\text{with} \quad  \boldsymbol{\delta} & = (\delta_{1}, ..., \delta_{N_{l}})^{T}, \quad \delta_{i} \sim \mathcal{N}(0, \sigma_{i}^{2}).
\end{align*}

We carry out our convergence analysis under following assumption.

\begin{asn}
\label{covariance}
\textit{The covariance matrix $\boldsymbol{\Sigma_{\sigma}}$ has the simple form} 
\begin{equation*}
\boldsymbol{\Sigma_{\sigma}} = \delta^{2}\cdot \mathrm{diag}(\sigma_{1}^{2}, ..., \sigma_{N_{l}}^{2}) =: \delta^{2} \cdot \boldsymbol{\Sigma},
\end{equation*}
\textit{where $\delta \geq 0$ is an arbitrary but fixed noise level and $\sigma_{1}$, ..., $\sigma_{N_{l}}$ are fixed.}
\end{asn}

Now instead of maximizing the posterior probability \eqref{post_prob} directly, we use the fact that
\begin{equation*}
\exp(-\textstyle{\frac{1}{2}}\|\boldsymbol{\Sigma}^{-\frac{1}{2}}(\boldsymbol{K}_{N}\boldsymbol{n} - (\boldsymbol{e}_{true} + \boldsymbol{\delta}))\|_{2}^{2} - \textstyle{\frac{1}{2}}\gamma\delta^{2}\|\boldsymbol{n}\|_{2}^{2})I_{S}(\boldsymbol{n})
\end{equation*}
has the same maximizer. To obtain the function above we scaled the argument of the exponential in \eqref{post_prob} with the noise level $\delta^{2}$. For simpler notation we redefine for all the following
\begin{equation*}
\boldsymbol{K} := \boldsymbol{\Sigma}^{-\frac{1}{2}}\boldsymbol{K}_{N}, \quad \boldsymbol{r} := \boldsymbol{\Sigma}^{-\frac{1}{2}}(\boldsymbol{e}_{true} + \boldsymbol{\delta}) \quad \text{and} \quad \alpha = \gamma\delta^{2}.
\end{equation*}
This means that we work with versions of $\boldsymbol{\Sigma_{\sigma}}^{-\frac{1}{2}}(\boldsymbol{e}_{true} + \boldsymbol{\delta})$ and $\boldsymbol{\Sigma_{\sigma}}^{-\frac{1}{2}}\boldsymbol{K}_{N}$ where the noise magnitude $\delta^{2}$ is scaled out.
 So instead of solving \eqref{quadprog}, we now solve
\begin{equation}
\min_{\boldsymbol{n} \in \mathbb{R}^{N}} \textstyle{\frac{1}{2}}\|\boldsymbol{K}\boldsymbol{n} - \boldsymbol{r}\|_{2}^{2} + \textstyle{\frac{1}{2}}\alpha\|\boldsymbol{n}\|_{2}^{2} \quad \text{s.t.} \quad \boldsymbol{Cn} \leq \boldsymbol{b}.
\label{quadprog2}
\end{equation}

We have to point out that the back-scaled parameter $\gamma = \alpha / \delta^{2}$ must be used for the statistical computations for the posterior probabilities in Section \ref{section_model_generation}. So we always compute the parameter $\alpha$ first from \eqref{quadprog2} and then obtain $\gamma$ from it. We can already see here that Bayesian model-selection computations are not feasible for very small noise levels $\delta$, since the parameter $\gamma$ diverges as $\delta$ tends to $0$. Another reason for skipping the model selection step for $\delta$ approaching $0$ is that the entries of the covariance matrix $\boldsymbol{\Sigma_{\sigma}}$ get closer to $0$ as well here, which causes problems in the statistical computations which will follow in Section \ref{section_model_generation}. We recommend to switch to the classical discrepany principle in this case.

Now we present the standard convergence rate for Tikhonov regularization.   

\begin{prp}
\label{convergence}
\textit{If the noise-free true solution $\boldsymbol{n}_{0}$ is an element of the feasible set of \eqref{quadprog}, then the regularized solutions $\boldsymbol{n}_{\alpha}$ of the noise-free problem satisfies}
\begin{equation}
\|\boldsymbol{K}(\boldsymbol{n}_{0} - \boldsymbol{n}_{\alpha})\|_{2} = \mathcal{O}(\alpha^{\frac{1}{2}}) 
\end{equation}
\textit{as $\alpha$ goes to $0$. Thus $\lim_{\alpha \to 0}\boldsymbol{n}_{\alpha} = \boldsymbol{n}_{0}$.}
\end{prp}

\begin{prf}
Rearranging \eqref{main_eq} for $\beta = 0$ gives the error representation
\begin{equation*}
\|\boldsymbol{K}(\boldsymbol{n}_{0} - \boldsymbol{n}_{\alpha})\|_{2}^{2} = \langle \boldsymbol{n}_{0} - \boldsymbol{n}_{\alpha}, \alpha \boldsymbol{n}_{\alpha} - \boldsymbol{C}^{T}(\boldsymbol{q}_{0} - \boldsymbol{q}_{\alpha}) \rangle.
\end{equation*}
Now since $ \langle \boldsymbol{n}_{0} - \boldsymbol{n}_{\alpha}, \boldsymbol{C}^{T}(\boldsymbol{q}_{0} - \boldsymbol{q}_{\alpha}) \rangle \geq 0$ and $\|\boldsymbol{n}_{\alpha}\|_{2} \leq \|\boldsymbol{n}_{0}\|_{2}$ hold, we can therefore estimate
\begin{equation*}
\|\boldsymbol{K}(\boldsymbol{n}_{0} - \boldsymbol{n}_{\alpha})\|_{2}^{2} \leq \alpha \|\boldsymbol{n}_{0}\|_{2}^{2}
\end{equation*}
which gives the first result. The second assertion was proved in \cite{Neub88a}.
\hfill $\square$
\end{prf}

\begin{prp}
\label{bounds}
\textit{Let $\boldsymbol{r}$ and $\boldsymbol{\tilde{r}}$ be two different data vectors for \eqref{quadprog2} and let $\boldsymbol{n}_{\alpha}$ and $\boldsymbol{\tilde{n}}_{\alpha}$ be the corresponding regularized solutions of \eqref{quadprog2} for the parameter $\alpha$. Then}
\begin{equation}
\|\boldsymbol{K}(\boldsymbol{n}_{\alpha} - \boldsymbol{\tilde{n}}_{\alpha})\|_{2} \leq \|\boldsymbol{r} - \boldsymbol{\tilde{r}}\|_{2}  \quad \text{and} \quad\|\boldsymbol{n}_{\alpha} - \boldsymbol{\tilde{n}}_{\alpha}\|_{2} \leq \frac{\|\boldsymbol{r} - \boldsymbol{\tilde{r}}\|_{2}}{\alpha^{\frac{1}{2}}}.
\end{equation}
\end{prp}

\begin{prf}
We give the proof from \cite{Neub88a}. The solutions $\boldsymbol{n}_{\alpha}$ and $\boldsymbol{\tilde{n}}_{\alpha}$ fulfill the variational inequalities
\begin{align*}
\langle \boldsymbol{K}^{T}\boldsymbol{K}\boldsymbol{n}_{\alpha} -  \boldsymbol{K}^{T}\boldsymbol{r} + \alpha \boldsymbol{n}_{\alpha}, \boldsymbol{\tilde{n}}_{\alpha} - \boldsymbol{n}_{\alpha} \rangle & \geq 0 \\
\text{and} \quad \langle \boldsymbol{K}^{T}\boldsymbol{K}\boldsymbol{\tilde{n}}_{\alpha} -  \boldsymbol{K}^{T}\boldsymbol{\tilde{r}} + \alpha \boldsymbol{\tilde{n}}_{\alpha}, \boldsymbol{n}_{\alpha} - \boldsymbol{\tilde{n}}_{\alpha} \rangle & \geq 0.
\end{align*}
Adding them gives
\begin{align*}
\|\boldsymbol{K}(\boldsymbol{\tilde{n}}_{\alpha} - \boldsymbol{n}_{\alpha})\|_{2}^{2} + \alpha \|\boldsymbol{\tilde{n}}_{\alpha} - \boldsymbol{n}_{\alpha}\|_{2}^{2} & \leq \langle \boldsymbol{\tilde{r}} - \boldsymbol{r}, \boldsymbol{K}(\boldsymbol{\tilde{n}}_{\alpha} - \boldsymbol{n}_{\alpha})\rangle \\
& \leq \|\boldsymbol{\tilde{r}} - \boldsymbol{r}\|_{2}\|\boldsymbol{K}(\boldsymbol{\tilde{n}}_{\alpha} - \boldsymbol{n}_{\alpha})\|_{2},
\end{align*}
and the desired results follow from the last inequality.
\hfill $\square$
\end{prf}

Finally we show under which conditions the regularized solutions $\boldsymbol{n}_{\alpha}^{\delta}$ of the noisy problem \eqref{quadprog2} converge to the true solution $\boldsymbol{n}_{0}$ of the noise-free problem for $\delta \rightarrow 0$. In preparation we note that for the weighted residual with noise level $\boldsymbol{\delta}$ 
\begin{align*}
\|\boldsymbol{\Sigma_{\sigma}}^{-\frac{1}{2}}\boldsymbol{\delta}\|_{2}^{2} \sim \chi^{2}(N_{l}) \quad \text{thus} \quad & \mathbb{E}\big(\|\boldsymbol{\Sigma}^{-\frac{1}{2}}\boldsymbol{\delta}\|_{2}^{2}\big) =  \delta^{2} \cdot \mathbb{E}\big(\|\boldsymbol{\Sigma_{\sigma}}^{-\frac{1}{2}}\boldsymbol{\delta}\|_{2}^{2}\big) = N_{l}\delta^{2} = \mathcal{O}(\delta^{2}), \\
\text{i.e.} \quad & \mathbb{E}\big(\|\boldsymbol{\Sigma}^{-\frac{1}{2}}\boldsymbol{\delta}\|_{2}\big) \leq \left(\mathbb{E}\big(\|\boldsymbol{\Sigma}^{-\frac{1}{2}}\boldsymbol{\delta}\|_{2}^{2}\big)\right)^{\frac{1}{2}} = \mathcal{O}(\delta).
\end{align*}
We set $\boldsymbol{r}_{true} := \boldsymbol{\Sigma}^{-\frac{1}{2}}\boldsymbol{e}_{true}$. Then for the expected value we have
\begin{equation*}
\mathbb{E}(\|\boldsymbol{r} - \boldsymbol{r}_{true}\|_{2}) = \mathbb{E}(\|\boldsymbol{\Sigma}^{-\frac{1}{2}}\boldsymbol{\delta}\|_{2}) = \mathcal{O}(\delta).
\end{equation*}

\begin{thm}
\label{conv_thm1}
\textit{If we have $\mathbb{E}(\|\boldsymbol{r} - \boldsymbol{r}_{true}\|_{2}) = \mathcal{O}(\delta)$ and $\alpha(\delta)$ has the properties
$\lim_{\delta \to 0} \: \alpha(\delta) = 0$ and $\lim_{\delta \to 0} \: \frac{\delta^{2}}{\alpha(\delta)} = 0,$ then $\lim_{\delta \to 0} \: \mathbb{E}\big(\|\boldsymbol{n}_{\alpha(\delta)}^{\delta} - \boldsymbol{n}_{0}\big\|_{2}) = 0$ holds.} 
\end{thm}

\begin{prf}
We have
\begin{equation*}
\mathbb{E}\big(\|\boldsymbol{n}_{\alpha(\delta)}^{\delta} - \boldsymbol{n}_{0}\|_{2}\big) \leq \mathbb{E}\big(\|\boldsymbol{n}_{\alpha(\delta)}^{\delta} - \boldsymbol{n}_{\alpha(\delta)}\|_{2}\big) + \mathbb{E}(\|\boldsymbol{n}_{\alpha(\delta)} - \boldsymbol{n}_{0}\|_{2}),
\end{equation*}
where $\boldsymbol{n}_{\alpha(\delta)}$ is the regularized solution for the noise-free data $\boldsymbol{r}_{true}$. Having $\mathbb{E}(\|\boldsymbol{r} - \boldsymbol{r}_{true}\|_{2}) = \mathcal{O}(\delta)$, we can further estimate using Proposition \ref{bounds}
\begin{equation*}
\mathbb{E}\big(\|\boldsymbol{n}_{\alpha(\delta)}^{\delta} - \boldsymbol{n}_{0}\|_{2}\big) \leq \mathbb{E}(\|\boldsymbol{n}_{\alpha(\delta)} - \boldsymbol{n}_{0}\|_{2}) + \frac{\mathcal{O}(\delta)}{\alpha^{\frac{1}{2}}}.
\end{equation*}
Then the result follows with Proposition \ref{convergence}.
\hfill $\square$
\end{prf}

\subsection{Model Generation under Nonnegativity Constraints}
\label{section_model_generation}

Suppose we have discretized our linear operator with a Galerkin collocation method on a set of $m$ different grids. Each grid has $N_{k}$ collocation points with $N_{1} < ... < N_{m}$ and we have computed a discrete approximation $\boldsymbol{K}_{k}$ to $K$ for each grid. The approximation $\boldsymbol{n}_{k}$ of the sought-after function $n$ lies in $\mathbb{R}^{N_{k}}$. For each grid we apply a Tikhonov prior with nonnegativity constraints on the observed model uncertainty such that we have according to the previously derived results a bijection between attainable residuals and regularization parameters.

Because $\delta_{1}, ..., \delta_{N_{l}}$ are normally distributed, it follows with $\boldsymbol{\Sigma_{\sigma}}^{-\frac{1}{2}}\boldsymbol{K}_{k}\boldsymbol{n}_{k} = \boldsymbol{e}_{true}$ that
\begin{align*}
& \|\boldsymbol{\Sigma_{\sigma}}^{-\frac{1}{2}}\left(\boldsymbol{K}_{k}\boldsymbol{n}_{k} - (\boldsymbol{e}_{true} + \boldsymbol{\delta})\right)\|_2^{2} \sim \chi^{2}(N_{l}), \\
\text{and thus} \quad \mathbb{E}\big(& \|\boldsymbol{\Sigma_{\sigma}}^{-\frac{1}{2}}\left(\boldsymbol{K}_{k}\boldsymbol{n}_{k} - (\boldsymbol{e}_{true} + \boldsymbol{\delta})\right)\|_2^{2}\big) = N_{l}, \quad \forall k \in \{1, ..., m\}.
\end{align*} 
In the literature on the discrepancy principle, e.g. in \cite{KS05}, the error estimate $N_{l}$ is multiplied with a factor $\tau$ near $1$ which is known as \textit{Morozov's safety parameter}. Now we interpret it here statistically as high-probability values of the observed distribution of the weighted residual. Of course we do not select just one single value for $\tau$, instead we select a grid of Morozov safety parameters $\tau_{1}, ..., \tau_{s}$. The following example of the $\chi^{2}(48)$ probability density functions illustrates this strategy:

\begin{figure}[h]
\centering
\includegraphics[width =1.0\textwidth]{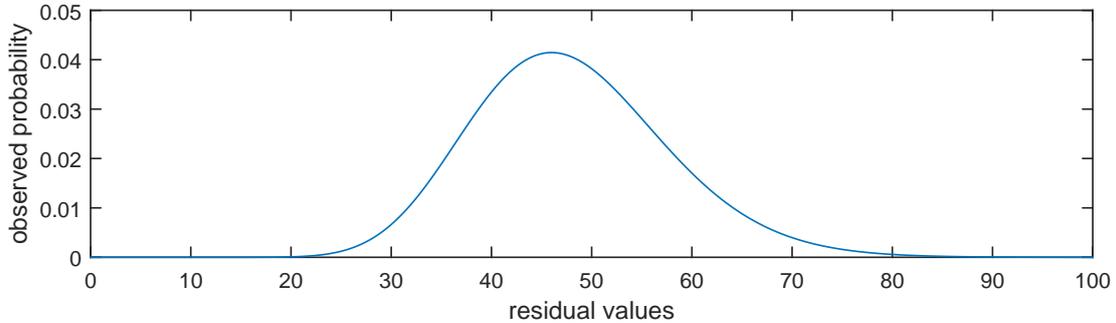}
\caption{$\chi^{2}(48)$ probability density function}
\label{chisqu48}
\end{figure}

\noindent It is indeed a unimodal distribution with residual values having a nonnegliglible probability ranging from $30$ to $70$. For $N_{l} = 48$ this corresponds to values of $\tau$ ranging from ca. $0.6$ to ca. $1.5$. Therefore proposing just a single residual value for the discrepancy principle ($1.1 N_{l}$ would be a common choice) excludes many probable reconstructions corresponding to other residual values, such that the posterior probability exploration is limited. Moreover the danger of under- or overregularization would be high.

As in the previous section we use the normalized version $\boldsymbol{\Sigma}$ of the covariance matrix $\boldsymbol{\Sigma_{\sigma}}$. This means that we try to fit the normalized residuals $\|\boldsymbol{\Sigma}^{-\frac{1}{2}}\left(\boldsymbol{K}_{k}\boldsymbol{n}_{k} - (\boldsymbol{e}_{true} + \boldsymbol{\delta})\right)\|_2^{2}$ to the values $\tau N_{l} \delta^{2}$ where the values for $\tau$ run through the grid of preselected Morozov safety factors. In practice the noise magnitude $\delta^{2}$ is taken as the biggest measurement sample mean and $\boldsymbol{\Sigma}$ is estimated from $\boldsymbol{\Sigma_{\sigma}}$ by normalizing it with the estimate for $\delta^{2}$.

For the following we set $\boldsymbol{e}_{real} := (\tilde{e}_{1}, ..., \tilde{e}_{N_{l}})^{T}$, hence this is the vector of the realizations of the random variables $e_{1}, ..., e_{N_{l}}$. With these preparations the model generation step proceeds as follows:

\newpage

\begin{algorithm}[h!]
\caption{Model Generation}
\label{model_generation}
\begin{algorithmic}[1]
\State $MaxDisc = 3$
\State $SolutionSets = \{\}$
\State $ApproxSets = \{\}$ 
\State $PriorSets = \{\}$
\State $TauSets = \{\}$
\State $DiscCntr = 0$
\State estimate $\sigma_{1}^{2}$, ..., $\sigma_{N_{l}}^{2}$ from the sample means approximating the standard deviations of $e_{1}$, ..., $e_{N_{l}}$. \vspace{0.05cm}  
%
%
\State $\delta^{2} := \mathrm{max}\big\{\sigma_{1}^{2}, ..., \sigma_{N_{l}}^{2}\big\}$
\State $\boldsymbol{\Sigma} := \delta^{-2} \cdot \mathrm{diag}\big(\sigma_{1}^{2}, ..., \sigma_{N_{l}}^{2}\big)$
     \For{$i = 1 \; \textbf{to} \; m$}
     \State $S_{i} = \{\}$
     \State $A_{i} = \{\}$
     \State $P_{i} = \{\}$
     \State $T_{i} = \{\}$
     \State $\boldsymbol{n}_{lsqnng} = \underset{\boldsymbol{n} \: \in \: \mathbb{R}^{N_{i}}}{\mathrm{argmin}} \; \textstyle{\frac{1}{2}}\|\boldsymbol{\Sigma}^{-\frac{1}{2}}(\boldsymbol{K}_{i}\boldsymbol{n} - \boldsymbol{e}_{real})\|_2^{2} \;\; \text{s.t.} \; \boldsymbol{n} \geq 0$ \vspace{0.1cm}
     \State $R_{lsqnng} = \|\boldsymbol{\Sigma}^{-\frac{1}{2}}(\boldsymbol{K}_{i}\boldsymbol{n}_{lsqnng} - \boldsymbol{e}_{real})\|_2^{2}$ \vspace{0.05cm}  
	  \For{$j = 1 \; \textbf{to} \; s$}	
	         	\If{$R_{lsqnng} < \tau_{j}N_{l}\delta^{2} \; \land \; \tau_{j}N_{l}\delta^{2} < \|\boldsymbol{\Sigma}^{-\frac{1}{2}}\boldsymbol{e}_{real}\|_2^{2}$} \vspace{0.1cm} \label{DiscPrincipleCheck} 
	  	\State $\text{compute }\gamma_{ij}\text{ such that}$ \vspace{0.1cm}	
	  	\State $\boldsymbol{n}_{trial} = \underset{\boldsymbol{n} \: \in \: \mathbb{R}^{N_{i}}}{\mathrm{argmin}} \; \textstyle{\frac{1}{2}}\|\boldsymbol{\Sigma}^{-\frac{1}{2}}(\boldsymbol{K}_{i}\boldsymbol{n} - \boldsymbol{e}_{real})\|_2^{2} + \textstyle{\frac{1}{2}}\gamma_{ij}\boldsymbol{n}^{T}\boldsymbol{R}_{i}\boldsymbol{n} \;\; \text{s.t.} \; \boldsymbol{n} \geq 0$ \vspace{0.05cm}
	  \State $\text{with } \|\boldsymbol{\Sigma}^{-\frac{1}{2}}(\boldsymbol{K}_{i}\boldsymbol{n}_{trial} - \boldsymbol{e}_{real})\|_2^{2} = \tau_{j}N_{l}\delta^{2}$	\vspace{0.05cm}  
		\EndIf
	  	\If{$\boldsymbol{n}_{trial}\text{ exists}$}
		\State $S_{i} = S_{i} \cup \{\boldsymbol{n}_{trial}\}$
		\State $A_{i} = A_{i} \cup \{\boldsymbol{K}_{i}\}$
		\State $P_{i} = P_{i} \cup \{\gamma_{ij}\boldsymbol{R}_{i}\}$
		\State $T_{i} = T_{i} \cup \{\tau_{j}\}$
		\EndIf
	  \EndFor	
	  \If{$S_{i}, \; A_{i}, \; P_{i} \text{ and } T_{i} \text{ not empty}$}
	  \State $SolutionSets = SolutionSets \cup \{S_{i}\}$
	  \State $ApproxSets = ApproxSets \cup \{A_{i}\}$
	  \State $PriorSets = PriorSets \cup \{P_{i}\}$
	  \State $TauSets = TauSets \cup \{T_{i}\}$
	  \State $DiscCntr = DiscCntr + 1$
	  \EndIf     
	  \If{$DiscCntr == MaxDisc$}
	  \State $\text{break}$
	  \EndIf
     \EndFor
\end{algorithmic}
\end{algorithm}

The outer loop runs through the discretization levels beginning with the coarsest one. This approach is in accordance with the principle of \textit{Occam's razor}, where among all possible explanations of a problem simpler ones are preferred over more complicated ones. Another motivation is \textit{regularization by discretization}, which means that the approximate problems for the operator inversion are for coarser discretizations less ill-conditioned than for finer discretizations. But by using the discrepancy principle we ensure that the models selected are not too coarse by demanding that the model has to fit the data, which means that the residuals may not be too big.

For each $i$-th discretization level in the outer loop, the inner loop runs through the preselected grid of Morozov safety factors, where for each factor $\tau_{j}$ the computation of a regularized solution $\boldsymbol{n}_{trial}$ with residual $\tau_{j} N_{l}$ is attempted. In line \ref{DiscPrincipleCheck} it is checked if the discrepancy principle is applicable. If it is possible to compute $\boldsymbol{n}_{trial}$, this reconstruction is stored in the container $S_{i}$ and the approximation $\boldsymbol{K}_{i}$ to $K$ in $A_{i}$. The prior information given by the regularization parameter $\gamma_{ij}$ and the regularization matrix $\frac{1}{2}\boldsymbol{R}_{i}$ are stored in $P_{i}$ and the residual parameter $\tau_{j}$ in $T_{i}$. These matrices will be used to compute the Bayesian posterior probabilities for the model selection in the next section.

If in the current discretization level the containers with reconstructions, operator approximation matrices, prior informations and residual parameters are not empty, they are be added to the containers $SolutionSets$, $ApproxSets$, $PriorSets$ and $TauSets$ respectively. Note that we have limited the maximal number of admissible discretization levels to three. On the one hand this is done to save computational effort, but on the other hand it turns out that the posterior probabilities get too similar and thus not clearly or reliably distinguishable when using too many finely discretized models.

\section{Model Selection}
\label{section_model_selection}

\subsection{Posterior Model Probabilities under Nonnegativity Constraints}

In this section we apply the Bayesian model selection framework as introduced in \cite{Chi01}. Since we assume that the data is given by independent Gaussian random variables, the observed model uncertainty is a multivariate Gaussian distribution. For any of the approximations $\boldsymbol{K}_{k}$ to the operator $K$ with $k \in \{1, ..., m\}$ it is given by 
\begin{equation}
p(\boldsymbol{e} | \boldsymbol{n}, N_{k}, \boldsymbol{K}_{k}) = (2\pi)^{-\frac{N_{l}}{2}}\big|\mathrm{det}(\boldsymbol{\Sigma_{\sigma}})\big|^{-\frac{1}{2}}\exp(-\textstyle{\frac{1}{2}}(\boldsymbol{K}_{k}\boldsymbol{n} - \boldsymbol{e})^{T}\boldsymbol{\Sigma}^{-1}_{\boldsymbol{\sigma}}(\boldsymbol{K}_{k}\boldsymbol{n} - \boldsymbol{e})).
\label{model_uncertainty}
\end{equation}
Here the vector $\boldsymbol{n} \in \mathbb{R}^{N_{k}}$ represents all possible reconstructions for the current discretization.

We know beforehand that our reconstruction must be nonnegative and that it is smooth. We put this prior knowledge into our reconstruction method by setting up the Bayesian conditional prior probability which is determined by
\begin{equation}
p(\boldsymbol{n} | N_{k}, \boldsymbol{K}_{k}, \boldsymbol{R}_{k}, \gamma_{kj}) = C_{kj}^{-1}\exp(-\textstyle{\frac{1}{2}}\gamma_{kj}\boldsymbol{n}^{T}\boldsymbol{R}_{k}\boldsymbol{n})I_{\geq 0}(\boldsymbol{n}),
\label{prior}
\end{equation}
where $I_{\geq 0}(\boldsymbol{n})$ is the indicator function of the first quadrant of $\mathbb{R}^{N_{k}}$, $\boldsymbol{R}_{k}$ is the regularization matrix and $\gamma_{kj}$ is the regularization parameter. All these quantities were computed and stored in the model generation procedure in the previous section.

If $\boldsymbol{R}_{k}$ is regular and positive definite, the normalizing constant 
\begin{equation}
C_{kj} = \int_{[0, \infty)^{N_{k}}}\exp(-\textstyle{\frac{1}{2}}\gamma_{kj}\boldsymbol{n}^{T}\boldsymbol{R}_{k}\boldsymbol{n}) d\boldsymbol{n}
\label{normalizing_constant}
\end{equation} 
is well-defined. For \textit{Tikhonov regularization}, where $\boldsymbol{R}_{k} = \boldsymbol{I}_{N_{k}}$ holds, we have a closed form expression for it, namely
\begin{equation*}
C_{kj} = \left(\frac{\pi}{2\gamma_{kj}}\right)^{\frac{N_{k}}{2}}.
\end{equation*}
 
In \textit{minimal first differences regularization} with zero boundary conditions the regularization matrix is given by
\begin{equation*}
\boldsymbol{R}_{k} = \boldsymbol{H}_{k}^{T}\boldsymbol{H}_{k} \quad \text{with} \quad
\boldsymbol{H}_{k} = 
\begin{pmatrix}
-1 & & & \\
1 & -1 &  &  \\
& \ddots & \ddots & \\
& & 1 & -1  \\
& & & 1  \\
\end{pmatrix}.
\end{equation*} 

For \textit{Twomey regularization} with eliminated zero boundary conditions we have
\begin{equation*}
\boldsymbol{R}_{k} = \boldsymbol{H}_{k}^{T}\boldsymbol{H}_{k} \quad \text{with} \quad
\boldsymbol{H}_{k} = 
\begin{pmatrix}
2 & -1 & & &\\
-1 & 2 & -1 & & \\
& \ddots & \ddots & \ddots & \\
& & -1 & 2 & -1 & \\
& & & -1 & 2 
\end{pmatrix}.
\end{equation*} 
Here $\boldsymbol{R}_{k}$ is a positive definite tridiagonal matrix. For the latter two regularization methods $C_{kj}$ must be computed numerically. 
 
Remember that the container $SolutionSets$ stores reconstructions from at most $3$ discretization levels. We let the index $i$ run through all discretization levels in $SolutionSets$ and the index $j$ through all residual parameters captured in the $i$-th level. Then with Bayes' rule the posterior model probabilities are
\begin{equation*}
p(N_{k}, \gamma_{kt} | \boldsymbol{e}) 
= \frac{p(\boldsymbol{e} | N_{k}, \gamma_{kt}) p(N_{k}, \gamma_{kt})}{\sum_{i}\sum_{j}p(\boldsymbol{e} | N_{i}, \gamma_{ij}) p(N_{i}, \gamma_{ij})}
\end{equation*}
where with \eqref{model_uncertainty}-\eqref{normalizing_constant} we have
\begin{equation}
\begin{split}
& p(\boldsymbol{e} | N_{i}, \gamma_{ij}) \\
= & \displaystyle{\int_{\mathbb{R}^{N_{i}}}} p(\boldsymbol{e}, \boldsymbol{n} | N_{i}, \gamma_{ij}) d\boldsymbol{n} \\
= & \displaystyle{\int_{\mathbb{R}^{N_{i}}}} p(\boldsymbol{e} | \boldsymbol{n}, N_{i}) p(\boldsymbol{n} | N_{i}, \gamma_{ij}) d\boldsymbol{n} \\
= & \displaystyle{\int_{[0, \infty)^{N_{i}}}B^{-1}C_{ij}^{-1}\exp(-\textstyle{\frac{1}{2}}\big\|\boldsymbol{\Sigma}^{-\frac{1}{2}}_{\boldsymbol{\sigma}}(\boldsymbol{K}_{i}\boldsymbol{n} - \boldsymbol{e})\big\|_2^{2} - \textstyle{\frac{1}{2}}\gamma_{ij}\boldsymbol{n}^{T}\boldsymbol{R}_{i}\boldsymbol{n})d\boldsymbol{n}}, \\
\end{split}
\label{marginal}
\end{equation}
where
\begin{align*}
B & = (2\pi)^{\frac{N_{l}}{2}}\big|\mathrm{det}(\boldsymbol{\Sigma_{\sigma}})\big|^{\frac{1}{2}} \\
& \\
\text{and} \quad C_{ij} & = \int_{[0, \infty)^{N_{i}}}\exp(-\textstyle{\frac{1}{2}}\gamma_{ij}\boldsymbol{n}^{T}\boldsymbol{R}_{i}\boldsymbol{n}) d\boldsymbol{n}.
\end{align*}
We assumed that the model matrix $\boldsymbol{K}_{i}$ and the regularization matrix $\boldsymbol{R}_{i}$ were implicitly given by each discretization level $N_{i}$, i.e.\ we actually have $p(\boldsymbol{e} | N_{i}, \gamma_{ij}) = p(\boldsymbol{e} | N_{i}, \boldsymbol{K}_{i}, \boldsymbol{R}_{i}, \gamma_{ij})$. For simplicity of notation these were omitted. Note that the prior model probabilities $p(N_{i}, \gamma_{ij})$ are still free.

For the computation of the above integrals of multivariate Gaussian densities over the first quadrant of each model space $\mathbb{R}^{N_{i}}$ we applied an effective pseudo-random integration method described in \cite{Gen92} which implements the routines presented in \cite{Nie72} and \cite{CP76}.

\subsection{Model Selection for Nonnegativity Constraints}
\label{subsetion_model_selection}

We now turn to the prior model probabilities $p(N_{i}, \gamma_{ij})$. As mentioned in the beginning of Section \ref{constrained_Tikhonov} we assume that a $\gamma_{min} > 0$ and a $\gamma_{max} < \infty$ exist which give a lower and an upper bound for the regularization parameters $\gamma$ in order to exclude improper or point-mass priors for the cases $\gamma = 0$ or $\gamma = \infty$. This assumption is independent of the discretization level. We further assume the discretization level to be independent and uniformly distributed.  Thus we are taking a \emph{noninformative} prior, and so the prior model probabilities cancel out and do not affect the posterior probabilities.

Now everything is prepared to perform the model selection. To compute integrals of the form 
\begin{equation*}
\int_{[0, \infty)^{N}} \exp(-\textstyle{\frac{1}{2}}(\boldsymbol{n}^{T}\boldsymbol{Hn} - 2\boldsymbol{n}^{T}\boldsymbol{v} + q)) d\boldsymbol{n}, 
\end{equation*}
where $N$ is the dimension of the square matrix $\boldsymbol{H}$, we apply the method from \cite{Gen92}. It actually can only evaluate intgrals of the form
\begin{equation*}
\frac{1}{\sqrt{\mathrm{det}(\boldsymbol{W})(2\pi)^{N}}}\int_{a_1}^{b_1} ... \int_{a_N}^{b_N}\exp(-\textstyle{\frac{1}{2}}\boldsymbol{n}^{T}\boldsymbol{W}^{-1}\boldsymbol{n}) d\boldsymbol{n}, 
\end{equation*} 
where the cases $a_i = - \infty$ and $b_i = \infty$ are allowed. So we have to perform a simple affine transformation using the Cholesky factorization $\boldsymbol{H} = \boldsymbol{U}^{T}\boldsymbol{U}$:
\begin{align*}
& \int_{[0, \infty)^{N}} \exp(-\textstyle{\frac{1}{2}}(\boldsymbol{n}^{T}\boldsymbol{Hn} - 2\boldsymbol{n}^{T}\boldsymbol{v} + q)) d\boldsymbol{n} \\
= \; & \left(\exp(-\textstyle{\frac{1}{2}}(q - \boldsymbol{v}^{T}\boldsymbol{H}^{-1}\boldsymbol{v}))\sqrt{\mathrm{det}(\boldsymbol{H}^{-1})(2\pi)^{N}}\right) \\
\cdot \; & \frac{1}{\sqrt{\mathrm{det}(\boldsymbol{H}^{-1})(2\pi)^{N}}}\int_{\big\{\boldsymbol{z} \in \mathbb{R}^{N} | \; \boldsymbol{z} \; \geq \; - \boldsymbol{H}^{-1}\boldsymbol{v}\big\}}\exp(-\textstyle{\frac{1}{2}}\boldsymbol{z}^{T}\boldsymbol{Hz}) d\boldsymbol{z}.
\end{align*}
The model selection algorithm is as follows.

\begin{algorithm}[H]
\caption{Model Selection}
\label{model_selection}
\begin{algorithmic}[1]
\State $\text{get } S_{1}, ..., S_{MaxDisc} \text{ from } SolutionSets$
\State $\text{get } A_{1}, ..., A_{MaxDisc} \text{ from } ApproxSets$
\State $\text{get } P_{1}, ..., P_{MaxDisc} \text{ from } PriorSets$
\State $\text{get } T_{1}, ..., T_{MaxDisc} \text{ from } TauSets$ \vspace{0.1cm}
\State $m_{1} = |S_{1}|, ..., m_{MaxDisc} = |S_{MaxDisc}|$ \vspace{0.15cm}
\State $m_{total} = \sum_{k=1}^{MaxDisc}m_{k}$ \vspace{0.1cm}
\State $B = (2\pi)^{\frac{N_{l}}{2}}\big|\mathrm{det}(\boldsymbol{\Sigma_{\sigma}})\big|^{\frac{1}{2}}$ \vspace{0.1cm}
\State $P_{post} = \{\}$ \vspace{0.1cm}
     \For{$i = 1 \; \textbf{to} \; MaxDisc$} \label{BeginSecondLoop}
          \For{$j = 1 \; \textbf{to} \; m_{i}$}
          \State $\boldsymbol{K}_{ij} = A_{i}(j)$ \vspace{0.1cm}
          \State $\boldsymbol{R}_{ij} = \frac{1}{\delta^{2}}P_{i}(j)$ \vspace{0.1cm} \label{RegMatrix}
%
%
          \State $C_{ij} = \displaystyle{\int_{[0, \infty)^{N_{i}}}}\exp(-\textstyle{\frac{1}{2}}\boldsymbol{n}^{T}\boldsymbol{R}_{ij}\boldsymbol{n}) d\boldsymbol{n}$ \vspace{0.1cm}
          \State $M_{ij} =  \displaystyle{\int_{[0, \infty)^{N_{i}}}}\exp(-\textstyle{\frac{1}{2}}\big\|\boldsymbol{\Sigma}^{-\frac{1}{2}}_{\boldsymbol{\sigma}}(\boldsymbol{K}_{ij}\boldsymbol{n} -  \boldsymbol{e}_{real})\big\|_2^{2} - \textstyle{\frac{1}{2}}\boldsymbol{n}^{2}\boldsymbol{R}_{ij}\boldsymbol{n}) d\boldsymbol{n}$ \vspace{0.1cm}
          \State $P_{post} = P_{post} \cup \big\{M_{ij}/(B \cdot C_{ij})\big\}$ \label{posterior}
          \EndFor
     \EndFor \label{EndSecondLoop} \vspace{0.1cm}
\State $SumP_{post} = \sum_{k=1}^{m_{total}}P_{post}(k)$ \label{BeginNormalize} \vspace{0.1cm}    
     \For{$i = 1 \; \textbf{to} \; MaxDisc$}
     \State $P_{post}(i) = P_{post}(i)/SumP_{post}$
     \EndFor \label{EndNormalize}
\State $S_{total} = S_{1} \cup ... \cup S_{MaxDisc}$
\State $\textbf{sort } S_{total}(1), ..., S_{total}(m_{total}) \textbf{ according to } P_{post}(1), ..., P_{post}(m_{total})$
\end{algorithmic}
\end{algorithm}

In the first lines of the model-selection algorithm the containers for computed reconstructions, operator approximation matrices, prior matrices and residual parameters are loaded for each examined discretization level. They store the results of the model-generation algorithm from Section \ref{section_model_generation}. In the case of too noisy or improper data it might happen that in the model generation step none of the models can fit the data. Then all containers are empty and the model selection algorithm has to be aborted. For simplicity we assume that the model-generation step was successful.

The double loop in lines \ref{BeginSecondLoop}-\ref{EndSecondLoop} performs the multidimensional integrations needed in \eqref{normalizing_constant} and \eqref{marginal}. In line \ref{posterior} these integrals are used for the unnormalized posterior probabilities $p(N_{i}, \gamma_{ij} | \boldsymbol{e})$ from \eqref{marginal}. Note that the prior model probabilities $p(N_{i}, \gamma_{ij})$ do not appear in the algorithm, since they are selected to be uniform and thus cancel out in the normalizing step performed in lines \ref{BeginNormalize} - \ref{EndNormalize}. At last all reconstructions are sorted according to their posterior probabilities.

We have to be careful not to forget to normalize the regularization matrices $P_{i}(j)$ with the estimated noise level $\delta^{2}$ as in line \ref{RegMatrix} because all statistical computations have to be carried out using the unnormalized covariance matrix $\boldsymbol{\Sigma_{\sigma}}$. For very small noise levels we recommend to skip the model selection step completely due to instabilities in the statistical computations mentioned above. It is sufficient to use only the coarsest model generated with the commonly used value $\tau = 1.1$ in this case.

\section{Numerical Results}

\subsection{Simulation of Aerosol Spectroscopy Measurements}
\label{simulated_measurements}

We applied our algorithm to a simplified version of problem \eqref{operator_equation}, where we assumed that we know the minimal and maximal particle radii $r_{min}$ and $r_{max}$. This led to the integral equation
\begin{equation}
\int_{r_{min}}^{r_{max}} k(r, l) n(r) dr = e(l).
\label{integral_equation}
\end{equation}
For the kernel function $k(r,l)$ from Mie theory we selected $\mathrm{H}_{2}\mathrm{O}$ as the material for the scattering particles and air for the medium.

In our simulations we assumed $r_{min} = 0.01 \; \mu$m and $r_{max} = 7.0 \; \mu$m. In practice the extinction function can only be measured for a finite number of light wavelengths $l_{1}, ..., l_{N_{l}}$. In our simulations we used the grid of $48$ wavelengths composed of $8$ linearly spaced wavelengths from $0.6 - 0.8 \; \mu$m, $8$ from $1.1 - 1.3 \; \mu$m, $8$ from $1.6 - 1.8 \; \mu$m, $16$ from $2.1 - 2.5 \; \mu$m and $8$ from $3.1 - 3.3 \; \mu$m. These five intervals were chosen to exclude wavelengths where light absorption by ambient water can occur which distorts the measured extinctions $e(l)$ heavily. That is, the selected wavelengths cover the so-called \textit{optical window} which is free from this unwanted physical effect.

We generated artificial extinction values $e(l_{i})$ for the selected $l_{1}, ..., l_{N_{l}}$ by solving the forward problem, which means inserting an original `true' size distributions $n(r)$ into the integral equation \eqref{integral_equation}. To avoid the \textit{inverse crime} we used a very fine grid with $10001$ points and the composite Simpson rule to compute the resulting integrals.

In each simulation run we generated a set of $300$ noisy extinctions from the artificial true extinction values by adding zero-mean Gaussian noise where the standard deviations were taken to be $30\%$ of the true extinction values $e(l_{i})$. This means that a vector $\boldsymbol{e}$ of noisy extinctions for each single measurement was modeled as
\begin{equation*}
(\boldsymbol{e})_{i} = e(l_{i}) + \delta_{i} \quad \text{ with } \; \delta_{i} \sim \mathcal{N}(0,(0.3 \cdot e(l_{i}))^{2}), \quad i = 1, ..., N_{l}.
\end{equation*}
We used the sample means and variances of these $300$ artificial noisy extinctions to do inferences about the simulated Gaussian noise.

For the discretization of \eqref{integral_equation} we used a Galerkin collocation method with linear basis functions on an integration grid with $N_{r} = 300$ equidistant points. We generated our model spaces by selecting collocation grids as near equidistant subgrids of the integration grid where the number of grid points $N_{col}$ ranged from $3$ (coarsest discretization level) to $50$ (finest discretization level). For the collocation grids we set up linearly spaced `pre-collocation grids' with $N_{col}$ points first and then performed a nearest-neighbor-fitting of their points to the integration grid, such that they became subgrids. Since we are considering size distributions which attain small values at the minimal and maximal radii, we assumed zero boundary conditions. This effectively reduced the number of unknowns $N$ in each model space from $N = 3, ..., 50$ to $N = 1, ..., 48$ and---more importantly---prevented the reconstructed size distributions from sheering out at the smallest radius value, which would have been a not reasonable behavior, physically speaking. It was important that the dimension $N$ of each model space never succeeded the number of measurements $N_{l} = 48$, such that the resulting regression problems were fully or overdetermined.

Let $r_{1}, ..., r_{N_{r}}$ denote the integration grid points. Let $\{ r_{1} = c_{1} < ... < c_{N_{col}} = r_{N_{r}} \} \subset \{ r_{1}, ..., r_{N_{r}} \}$ be a collocation grid. The \textit{triangular basis funktions} $b_{k}(r)$, $k = 1, ..., N_{col}$ are the piecewise linear functions on the intervals $[c_{1}, c_{2}]$, ..., $[c_{N_{col} - 1}, c_{N_{col}}]$ which fulfill
\begin{equation*}
b_{k}(c_{j}) = \delta_{kj}, \text{ for }  j = 1, ..., N_{col}.
\end{equation*}
We approximated the sought-after function $n(r)$ with the linear combination
\begin{equation}
n(r) = \sum_{k = 1}^{N_{col}}n_{k}b_{k}(r),
\label{discretization}
\end{equation}
where the weights $n_{2}, ..., n_{N_{col} - 1} \in \mathbb{R}$ are free variables and $n_{1} = n_{N_{col}} = 0$ holds because of the zero boundary conditions.

Inserting \eqref{discretization} into \eqref{integral_equation} yields the linear system of equations for the unknown weights
\begin{equation}
\sum_{k = 1}^{N_{col}} n_{k} \int_{r_{min}}^{r_{max}} k(r, l_{i}) b_{k}(r) dr = e(l_{i}), \quad  i = 1, ..., N_{l}.
\end{equation}
We applied the composite trapezoidal rule with the integration grid $r_{1}, ..., r_{N{r}}$ on the integrals defining the coefficients in above linear system. The resulting coefficient matrix is the matrix $\boldsymbol{K}_{N}$ from Section \ref{section_model_generation} which approximates the integral operator from the left-hand side of \eqref{integral_equation}.

\subsection{Numerical Study}
\label{num_study}

We performed a numerical study for our reconstruction algorithm with model size distributions characterized by a low number of parameters. We varied the parameters in domains giving physically reasonable size distributions and generated noise in the same order of magnitude as observed in real experimental FASP measurements. Therefore the numerical results should give good estimates of the quality of the reconstructions compared to real size distributions. In the same simulation runs we compared our algorithm with existing reconstruction methods.

\subsubsection{Applied Methods}

For all inversion methods applied in our numerical study we selected for the priors Tikhonov, minimal first differences and Phillips-Twomey regularization from Section \ref{section_model_selection}.

In our inversion method we set the Morozov safety factor grid to
\begin{equation*}
\tau_{1} = 0.6, \tau_{2} = 0.7, ..., \tau_{12} = 1.7. 
\end{equation*}
We refer to this as the \textit{constrained method} in the following.

To see that the constraints in the constrained method are worth the computational effort, we compared it with its counterpart without constraints, which we call the \textit{unconstrained method}. It performs the same model generation step based on the discrepancy principle with the same Morozov safety factors grid, but the constraints in \eqref{quadprog} were dropped. The computations for the model selection are much easier here, since the integrals of the multivariate Gaussian distributions over the parameter spaces can be evaluated analytically.   

By reducing the grid of Morozov safety factors in the constrained method simply to the classical value $\tau = 1.1$ we obtained another method participating in our numerical study. We call it the \textit{Morozov method}. The comparison with it shows whether the grid of Morozov safety factors is justified or not. 

We also implemented a classical model-selection method for the unconstrained problem which is independent of the prior. Here we compared the three coarsest models where the discrepancy principle was applicable with the \textit{Bayesian Information Criterion} (BIC), which was first introduced in \cite{Sch78}. The model with the lowest BIC-value
\begin{equation*}
-2\Big(-\textstyle{\frac{1}{2}}N_{l}\log(2\pi) - \textstyle{\frac{1}{2}}\log(\mathrm{det}(\boldsymbol{\Sigma_{\sigma}})) - \textstyle{\frac{1}{2}}\|\boldsymbol{\Sigma}^{-\frac{1}{2}}_{\boldsymbol{\sigma}}(\boldsymbol{K}_{N}\boldsymbol{n}_{ml} - \boldsymbol{e})\|_{2}^2\Big) + N\log(N_{l}),
\end{equation*}
where $\boldsymbol{n}_{ml}$ is the unconstrained maximum-likelihood solution, is selected here. We call this method the \textit{BIC method}.

\subsubsection{Model Size Distributions}
\label{model_size_distributions}

We generated the simulated measurement data vectors $\boldsymbol{e}_{true}$ by inserting one of the following three model size distributions adopted from \cite{RW08} into our integral equation \eqref{integral_equation}:

\begin{enumerate}
\item \textit{log-normal distribution}
\begin{equation}
n(r) = \frac{A}{\sqrt{2\pi}\sigma r} \exp \bigg( - \frac{1}{2\sigma^{2}} \big(\log (r) - \log (\mu)\big)^{2} \bigg)
\label{LogNormal}
\end{equation}
with amplitude $A$, standard deviation $\sigma$, and mean $\mu$.
\item \textit{Rosin-Rammler-Sperling-Bennet} (RRSB) \textit{distribution}
\begin{equation}
n(r) = \frac{AN}{\nu}\left(\frac{r}{\nu}\right)^{N - 1} \exp \bigg( - \left(\frac{r}{\nu}\right)^{N} \bigg)
\label{RRSB}
\end{equation}
with amplitude $A$, exponent $N$, and mean $\nu$.
\item \textit{Hedrih distribution}
\begin{equation}
n(r) = \frac{128 A r^{3}}{3 \eta^{4}} \exp \bigg( -\frac{4 r}{\eta}\bigg)
\label{Hedrih}
\end{equation}
with amplitude $A$ and mean $\mu$.
\end{enumerate}
For each simulated size distribution we set the amplitude to $A = 10^{4}$. We choose the remaining parameters so that the relation
\begin{equation}
n(r_{max}) \leq Tol
\label{tail_value}
\end{equation}
with $r_{max} = 7.0  \;\mu$m and $Tol = 10$ was satisfied.
This is to be consistent with the assumption, that we can neglect the tails of the distributions and truncate them at the maximal radius $r_{max}$. Furthermore we assumed the modal value of the log-normal and RRSB distributions to be greater or equal to $1.0 \; \mu m$ in order to exclude too peaked distributions.
For each of the above three model size distributions we looped in our simulations through a set of $100$ possible parameters satisfying \eqref{tail_value}.

For the log-normal distributions we first selected for the mean $\sigma$ a linearly spaced grid with ten points ranging from $0.2$ to $0.5$, i.e.\ $\sigma_{k} = 0.2 + 0.3\frac{k - 1}{9}$, $k = 1, ..., 10$. Then we saw after a lengthy calculation that \eqref{tail_value} is equivalent to
\begin{equation*}
r_{max}\exp \bigg( - \bigg(- 2 \sigma^{2} \log \bigg( \frac{\sqrt{2 \pi} r_{max}\sigma Tol}{A} \bigg)\bigg)^{\frac{1}{2}}\bigg) \geq \mu.
\end{equation*}
The modal value of the log-normal distribution is $r_{mod} = \exp\big(\log(\mu) - \sigma^{2}\big)$, so $r_{mod} \geq 1.0$ is equivalent to $\mu \geq 1.0 \exp\big(\sigma^{2}\big)$. Using the last two inequalities we selected
\begin{align*}
\mu_{kj} & = 1.0 \exp\big(\sigma_{k}^{2}\big) + \frac{j - 1}{9}\Big( v_{k}
 - 1.0 \exp\big(\sigma_{k}^{2}\big)\Big) \\ 
\text{with} \quad v_{k} & =  r_{max}\exp \bigg( - \bigg(- 2 \sigma_{k}^{2} \log \bigg( \frac{\sqrt{2 \pi} r_{max}\sigma_{k} Tol}{A} \bigg)\bigg)^{\frac{1}{2}}\bigg), \\
k & = 1, ..., 10, \; j = 1, ..., 10.
\end{align*}
These are the $100$ parameters used for the log-normal distributions.

For the RRSB distributions we took for the exponents $N$ the integer values $N_{k} = k + 2$, $k = 1, ..., 10$. We computed the auxiliary variables $p_{k}$ as the real-valued solutions of the equations
\begin{equation*}
p_{k}\exp (- p_{k}) = \frac{r_{max} Tol}{A N_{k}}
\end{equation*}
being greater than one. With some algebra one can see that 
\begin{equation*}
\nu \leq r_{max}\cdot p_{k}^{-\frac{1}{N_{k}}}
\end{equation*}
is then equivalent to \eqref{tail_value} for the RRSB distribution. The modal value of the RRSB distribution is $r_{mod} = \nu \left(\frac{N - 1}{N}\right)^{\frac{1}{N}}$, therefore $r_{mod} \geq 1.0$ is equivalent to $\nu \geq 1.0 \cdot \left(\frac{N - 1}{N}\right)^{-\frac{1}{N}}$. Using the last two inequalities we selected
\begin{align*}
\nu_{kj} & = 1.0 \cdot \left(\frac{N_{k} - 1}{N_{k}}\right)^{-\frac{1}{N_{k}}} + \frac{j - 1}{9}\left(  r_{max}\cdot p_{k}^{-\frac{1}{N_{k}}} - 1.0 \cdot \left(\frac{N_{k} - 1}{N_{k}}\right)^{-\frac{1}{N_{k}}}\right), \\
k & = 1, ..., 10, \; j = 1, ..., 10.
\end{align*}
Thus we have $100$ parameters for the RRSB distributions.

For the Hedrih distribution we found that \eqref{tail_value} is equivalent to $\eta \leq \eta_{max}$ with $\eta_{max} \approx 2.0566$. Thus we took for $\eta$ the values
\begin{equation*}
\eta_{k} = 0.8 + \frac{k - 1}{99}\left(\eta_{max} - 0.8\right), \quad k = 1, ..., 100.
\end{equation*}

For each of the three size distribution classes we simulated ten artificial noisy measurement-data vetors $\boldsymbol{e}$ as described in Section \ref{simulated_measurements} for each of the corresponding $100$ parameters. This resulted in total in $1000$ single simulated FASP experiments for one model size distribution class.

For every inversion we computed the $L^{2}$-error of the obtained reconstruction relative to the original size distribution and measured the total run time needed for the inversion. The computations were performed on a notebook with a $2.27$ GHz CPU and $3.87$ GB accessible primary memory.

\subsubsection{Average $L^{2}$-Errors}

\begin{table}[h!]
\centering
\begin{tabular}{||m{2.5cm}||m{2.5cm}|m{2.5cm}|m{2.5cm}||}
\hline
\multicolumn{4}{||c||}{\multirow{3}{*}{}} \tabularnewline[-2.5ex]
\multicolumn{4}{||c||}{Log-Normal Distribution} \tabularnewline
\multicolumn{4}{||c||}{} \tabularnewline[-2.5ex]
\hline
\multirow{4}{*}{\vspace{0.2cm}method} & \multicolumn{3}{c||}{} \tabularnewline[-2.5ex]
 & \multicolumn{3}{c||}{average $L^{2}$-errors (\%)} \tabularnewline
& \multicolumn{3}{c||}{} \tabularnewline[-2.5ex]
\cline{2-4}
& \centering \multirow{2}{*}{Tikhonov} & \centering min. first & \centering \multirow{2}{*}{Twomey} \tabularnewline
& & \centering fin. diff. & \tabularnewline
\hline
\hline
constrained & \centering 21.3917 & \centering 21.8413 & \centering 23.6893 \tabularnewline
\hline
Morozov & \centering 25.9880 & \centering 26.2520 & \centering 27.5610 \tabularnewline
\hline
unconstrained & \centering 31.1588 & \centering 34.1079 & \centering 37.7806 \tabularnewline
\hline
BIC & \centering 46.4902 & \centering 48.3897 & \centering 50.7711 \tabularnewline
\hline
\end{tabular}
\end{table}

\begin{table}[h!]
\centering
\begin{tabular}{||m{2.5cm}||m{2.5cm}|m{2.5cm}|m{2.5cm}||}
\hline
\multicolumn{4}{||c||}{\multirow{3}{*}{}} \tabularnewline[-2.5ex]
\multicolumn{4}{||c||}{RRSB Distribution} \tabularnewline
\multicolumn{4}{||c||}{} \tabularnewline[-2.5ex]
\hline
\multirow{4}{*}{\vspace{0.2cm}method} & \multicolumn{3}{c||}{} \tabularnewline[-2.5ex]
 & \multicolumn{3}{c||}{average $L^{2}$-errors (\%)} \tabularnewline
& \multicolumn{3}{c||}{} \tabularnewline[-2.5ex]
\cline{2-4}
& \centering \multirow{2}{*}{Tikhonov} & \centering min. first & \centering \multirow{2}{*}{Twomey} \tabularnewline
& & \centering fin. diff. & \tabularnewline
\hline
\hline
constrained & \centering 18.6192 & \centering 17.8924 & \centering 17.6709 \tabularnewline
\hline
Morozov & \centering 23.8631 & \centering 23.2843 & \centering 23.4669 \tabularnewline
\hline
unconstrained & \centering 29.1867 & \centering 33.2711 & \centering 37.6596 \tabularnewline
\hline
BIC & \centering 74.1802 & \centering 76.7546 & \centering 80.4552 \tabularnewline
\hline
\end{tabular}
\end{table}

\begin{table}[h!]
\centering
\begin{tabular}{||m{2.5cm}||m{2.5cm}|m{2.5cm}|m{2.5cm}||}
\hline
\multicolumn{4}{||c||}{\multirow{3}{*}{}} \tabularnewline[-2.5ex]
\multicolumn{4}{||c||}{Hedrih Distribution} \tabularnewline
\multicolumn{4}{||c||}{} \tabularnewline[-2.5ex]
\hline
\multirow{4}{*}{\vspace{0.2cm}method} & \multicolumn{3}{c||}{} \tabularnewline[-2.5ex]
 & \multicolumn{3}{c||}{average $L^{2}$-errors (\%)} \tabularnewline
& \multicolumn{3}{c||}{} \tabularnewline[-2.5ex]
\cline{2-4}
& \centering \multirow{2}{*}{Tikhonov} & \centering min. first & \centering \multirow{2}{*}{Twomey} \tabularnewline
& & \centering fin. diff. & \tabularnewline
\hline
\hline
constrained & \centering 14.3414 & \centering 13.2150 & \centering 12.8981 \tabularnewline
\hline
Morozov & \centering 25.6963 & \centering 24.3244 & \centering 23.5999 \tabularnewline
\hline
unconstrained & \centering 36.7919 & \centering 36.6560 & \centering 36.7877 \tabularnewline
\hline
BIC & \centering 42.7048 & \centering 41.8407 & \centering 41.5387 \tabularnewline
\hline
\end{tabular}
\end{table}

\newpage

\subsubsection{Average Run Times}
 
\begin{table}[h!]
\centering
\begin{tabular}{||m{2.5cm}||m{2.5cm}|m{2.5cm}|m{2.5cm}||}
\hline
\multicolumn{4}{||c||}{\multirow{3}{*}{}} \tabularnewline[-2.5ex]
\multicolumn{4}{||c||}{Log-Normal Distribution} \tabularnewline
\multicolumn{4}{||c||}{} \tabularnewline[-2.5ex]
\hline
\multirow{4}{*}{\vspace{0.2cm}method} & \multicolumn{3}{c||}{} \tabularnewline[-2.5ex]
 & \multicolumn{3}{c||}{average run times (s)} \tabularnewline
& \multicolumn{3}{c||}{} \tabularnewline[-2.5ex]
\cline{2-4}
& \centering \multirow{2}{*}{Tikhonov} & \centering min. first & \centering \multirow{2}{*}{Twomey} \tabularnewline
& & \centering fin. diff. & \tabularnewline
\hline
\hline
constrained & \centering 1.3525 & \centering 1.3987 & \centering 1.4068 \tabularnewline
\hline
Morozov & \centering 0.3136 & \centering 0.3284 & \centering 0.3325 \tabularnewline
\hline
unconstrained & \centering 0.1581 & \centering 0.1570 & \centering 0.1570 \tabularnewline
\hline
BIC & \centering 0.0324 & \centering 0.0324 & \centering 0.0327 \tabularnewline
\hline
\end{tabular}
\end{table}

\begin{table}[h!]
\centering
\begin{tabular}{||m{2.5cm}||m{2.5cm}|m{2.5cm}|m{2.5cm}||}
\hline
\multicolumn{4}{||c||}{\multirow{3}{*}{}} \tabularnewline[-2.5ex]
\multicolumn{4}{||c||}{RRSB Distribution} \tabularnewline
\multicolumn{4}{||c||}{} \tabularnewline[-2.5ex]
\hline
\multirow{4}{*}{\vspace{0.2cm}method} & \multicolumn{3}{c||}{} \tabularnewline[-2.5ex]
 & \multicolumn{3}{c||}{average run times (s)} \tabularnewline
& \multicolumn{3}{c||}{} \tabularnewline[-2.5ex]
\cline{2-4}
& \centering \multirow{2}{*}{Tikhonov} & \centering min. first & \centering \multirow{2}{*}{Twomey} \tabularnewline
& & \centering fin. diff. & \tabularnewline
\hline
\hline
constrained & \centering 1.7258 & \centering 1.8086 & \centering 1.8048 \tabularnewline
\hline
Morozov & \centering 0.3863 & \centering 0.4108 & \centering 0.4093 \tabularnewline
\hline
unconstrained & \centering 0.1523 & \centering 0.1516 & \centering 0.1512 \tabularnewline
\hline
BIC & \centering 0.0357 & \centering 0.0357 & \centering 0.0360 \tabularnewline
\hline
\end{tabular}
\end{table}

\begin{table}[h!]
\centering
\begin{tabular}{||m{2.5cm}||m{2.5cm}|m{2.5cm}|m{2.5cm}||}
\hline
\multicolumn{4}{||c||}{\multirow{3}{*}{}} \tabularnewline[-2.5ex]
\multicolumn{4}{||c||}{Hedrih Distribution} \tabularnewline
\multicolumn{4}{||c||}{} \tabularnewline[-2.5ex]
\hline
\multirow{4}{*}{\vspace{0.2cm}method} & \multicolumn{3}{c||}{} \tabularnewline[-2.5ex]
 & \multicolumn{3}{c||}{average run times (s)} \tabularnewline
& \multicolumn{3}{c||}{} \tabularnewline[-2.5ex]
\cline{2-4}
& \centering \multirow{2}{*}{Tikhonov} & \centering min. first & \centering \multirow{2}{*}{Twomey} \tabularnewline
& & \centering fin. diff. & \tabularnewline
\hline
\hline
constrained & \centering 1.2522 & \centering 1.2732 & \centering 1.2872 \tabularnewline
\hline
Morozov & \centering 0.2554 & \centering 0.2640 & \centering 0.2670 \tabularnewline
\hline
unconstrained & \centering 0.1572 & \centering 0.1565 & \centering 0.1561 \tabularnewline
\hline
BIC & \centering 0.0315 & \centering 0.0315 & \centering 0.0314 \tabularnewline
\hline
\end{tabular}
\end{table}

\newpage

\subsection{Average Model Space Dimensions}

\begin{table}[h!]
\centering
\begin{tabular}{||m{2.5cm}||m{2.5cm}|m{2.5cm}|m{2.5cm}||}
\hline
\multicolumn{4}{||c||}{\multirow{3}{*}{}} \tabularnewline[-2.5ex]
\multicolumn{4}{||c||}{Log-Normal Distribution} \tabularnewline
\multicolumn{4}{||c||}{} \tabularnewline[-2.5ex]
\hline
\multirow{4}{*}{\vspace{0.2cm}method} & \multicolumn{3}{c||}{} \tabularnewline[-2.5ex]
 & \multicolumn{3}{c||}{average model space dimensions} \tabularnewline
& \multicolumn{3}{c||}{} \tabularnewline[-2.5ex]
\cline{2-4}
& \centering \multirow{2}{*}{Tikhonov} & \centering min. first & \centering \multirow{2}{*}{Twomey} \tabularnewline
& & \centering fin. diff. & \tabularnewline
\hline
\hline
constrained & \centering 6.7570 & \centering 7.0890 & \centering 7.3670 \tabularnewline
\hline
Morozov & \centering 8.6010 & \centering 9.2220 & \centering 9.4540 \tabularnewline
\hline
unconstrained & \centering 5.2670 & \centering 5.3580 & \centering 5.4070 \tabularnewline
\hline
BIC & \centering 7.2310 & \centering 7.2310 & \centering 7.2310 \tabularnewline
\hline
\end{tabular}
\end{table}

\begin{table}[h!]
\centering
\begin{tabular}{||m{2.5cm}||m{2.5cm}|m{2.5cm}|m{2.5cm}||}
\hline
\multicolumn{4}{||c||}{\multirow{3}{*}{}} \tabularnewline[-2.5ex]
\multicolumn{4}{||c||}{RRSB Distribution} \tabularnewline
\multicolumn{4}{||c||}{} \tabularnewline[-2.5ex]
\hline
\multirow{4}{*}{\vspace{0.2cm}method} & \multicolumn{3}{c||}{} \tabularnewline[-2.5ex]
 & \multicolumn{3}{c||}{average model space dimensions} \tabularnewline
& \multicolumn{3}{c||}{} \tabularnewline[-2.5ex]
\cline{2-4}
& \centering \multirow{2}{*}{Tikhonov} & \centering min. first & \centering \multirow{2}{*}{Twomey} \tabularnewline
& & \centering fin. diff. & \tabularnewline
\hline
\hline
constrained & \centering 10.4220 & \centering 10.9830 & \centering 10.7100 \tabularnewline
\hline
Morozov & \centering 11.6910 & \centering 12.4760 & \centering 12.1860 \tabularnewline
\hline
unconstrained & \centering 8.7370 & \centering 9.5090 & \centering 9.2930 \tabularnewline
\hline
BIC & \centering 8.6840 & \centering 8.6840 & \centering 8.6840 \tabularnewline
\hline
\end{tabular}
\end{table}

\begin{table}[h!]
\centering
\begin{tabular}{||m{2.5cm}||m{2.5cm}|m{2.5cm}|m{2.5cm}||}
\hline
\multicolumn{4}{||c||}{\multirow{3}{*}{}} \tabularnewline[-2.5ex]
\multicolumn{4}{||c||}{Hedrih Distribution} \tabularnewline
\multicolumn{4}{||c||}{} \tabularnewline[-2.5ex]
\hline
\multirow{4}{*}{\vspace{0.2cm}method} & \multicolumn{3}{c||}{} \tabularnewline[-2.5ex]
 & \multicolumn{3}{c||}{average model space dimensions} \tabularnewline
& \multicolumn{3}{c||}{} \tabularnewline[-2.5ex]
\cline{2-4}
& \centering \multirow{2}{*}{Tikhonov} & \centering min. first & \centering \multirow{2}{*}{Twomey} \tabularnewline
& & \centering fin. diff. & \tabularnewline
\hline
\hline
constrained & \centering 6.0970 & \centering 6.3850 & \centering 6.6390 \tabularnewline
\hline
Morozov & \centering 7.5220 & \centering 7.9560 & \centering 8.1180 \tabularnewline
\hline
unconstrained & \centering 4.7820 & \centering 4.7980 & \centering 4.8200 \tabularnewline
\hline
BIC & \centering 7.4290 & \centering 7.4290 & \centering 7.4290 \tabularnewline
\hline
\end{tabular}
\end{table}

\subsection{Extreme Cases}

If the relative error of the reconstruction (compared with the original size distribution) is equal or even greater than $100$ percent, we regard the inversion as failed. Note that the inversion methods returned $\boldsymbol{n} \equiv 0$ by default if none of the kernel matrices in any of the model spaces would yield a reconstruction. Now we list how many times the inversion methods failed in our test runs. To see how trustworthy the results are we present the worst case $L^{2}$ errors as well. Finally we display the worst case run times.

\newpage

\subsubsection{Reconstruction Failures}

\begin{table}[h!]
\centering
\begin{tabular}{||m{2.5cm}||m{2.5cm}|m{2.5cm}|m{2.5cm}||}
\hline
\multicolumn{4}{||c||}{\multirow{3}{*}{}} \tabularnewline[-2.5ex]
\multicolumn{4}{||c||}{Log-Normal Distribution} \tabularnewline
\multicolumn{4}{||c||}{} \tabularnewline[-2.5ex]
\hline
\multirow{4}{*}{\vspace{0.2cm}method} & \multicolumn{3}{c||}{} \tabularnewline[-2.5ex]
 & \multicolumn{3}{c||}{number of $L^{2}$-errors $\geq 100$ \% (out of 1000)} \tabularnewline
& \multicolumn{3}{c||}{} \tabularnewline[-2.5ex]
\cline{2-4}
& \centering \multirow{2}{*}{Tikhonov} & \centering min. first & \centering \multirow{2}{*}{Twomey} \tabularnewline
& & \centering fin. diff. & \tabularnewline
\hline
\hline
constrained & \centering  0 & \centering  0 & \centering 0 \tabularnewline
\hline
Morozov & \centering 27 & \centering  27 & \centering  27 \tabularnewline
\hline
unconstrained & \centering 0 & \centering 0 & \centering 0 \tabularnewline
\hline
BIC & \centering 20 & \centering  20 & \centering 20 \tabularnewline
\hline
\end{tabular}
\end{table}

\begin{table}[h!]
\centering
\begin{tabular}{||m{2.5cm}||m{2.5cm}|m{2.5cm}|m{2.5cm}||}
\hline
\multicolumn{4}{||c||}{\multirow{3}{*}{}} \tabularnewline[-2.5ex]
\multicolumn{4}{||c||}{RRSB Distribution} \tabularnewline
\multicolumn{4}{||c||}{} \tabularnewline[-2.5ex]
\hline
\multirow{4}{*}{\vspace{0.2cm}method} & \multicolumn{3}{c||}{} \tabularnewline[-2.5ex]
 & \multicolumn{3}{c||}{number of $L^{2}$-errors $\geq 100$ \% (out of 1000)} \tabularnewline
& \multicolumn{3}{c||}{} \tabularnewline[-2.5ex]
\cline{2-4}
& \centering \multirow{2}{*}{Tikhonov} & \centering min. first & \centering \multirow{2}{*}{Twomey} \tabularnewline
& & \centering fin. diff. & \tabularnewline
\hline
\hline
constrained & \centering 2 & \centering 1 & \centering 1 \tabularnewline
\hline
Morozov & \centering 60 & \centering 60 & \centering 63 \tabularnewline
\hline
unconstrained & \centering 0 & \centering 0 & \centering 0 \tabularnewline
\hline
BIC & \centering 42 & \centering 47 & \centering 47 \tabularnewline
\hline
\end{tabular}
\end{table}

\begin{table}[h!]
\centering
\begin{tabular}{||m{2.5cm}||m{2.5cm}|m{2.5cm}|m{2.5cm}||}
\hline
\multicolumn{4}{||c||}{\multirow{3}{*}{}} \tabularnewline[-2.5ex]
\multicolumn{4}{||c||}{Hedrih Distribution} \tabularnewline
\multicolumn{4}{||c||}{} \tabularnewline[-2.5ex]
\hline
\multirow{4}{*}{\vspace{0.2cm}method} & \multicolumn{3}{c||}{} \tabularnewline[-2.5ex]
 & \multicolumn{3}{c||}{number of $L^{2}$-errors $\geq 100$ \% (out of 1000)} \tabularnewline
& \multicolumn{3}{c||}{} \tabularnewline[-2.5ex]
\cline{2-4}
& \centering \multirow{2}{*}{Tikhonov} & \centering min. first & \centering \multirow{2}{*}{Twomey} \tabularnewline
& & \centering fin. diff. & \tabularnewline
\hline
\hline
constrained & \centering 0 & \centering 0 & \centering 0 \tabularnewline
\hline
Morozov & \centering 33 & \centering 33 & \centering 33 \tabularnewline
\hline
unconstrained & \centering 0 & \centering 0 & \centering 0 \tabularnewline
\hline
BIC & \centering 12 & \centering 12 & \centering 12 \tabularnewline
\hline
\end{tabular}
\end{table}

\subsubsection{Worst Case Reconstruction Errors}

\begin{table}[h!]
\centering
\begin{tabular}{||m{2.5cm}||m{2.5cm}|m{2.5cm}|m{2.5cm}||}
\hline
\multicolumn{4}{||c||}{\multirow{3}{*}{}} \tabularnewline[-2.5ex]
\multicolumn{4}{||c||}{Log-Normal Distribution} \tabularnewline
\multicolumn{4}{||c||}{} \tabularnewline[-2.5ex]
\hline
\multirow{4}{*}{\vspace{0.2cm}method} & \multicolumn{3}{c||}{} \tabularnewline[-2.5ex]
 & \multicolumn{3}{c||}{worst case $L^{2}$-errors (\%)} \tabularnewline
& \multicolumn{3}{c||}{} \tabularnewline[-2.5ex]
\cline{2-4}
& \centering \multirow{2}{*}{Tikhonov} & \centering min. first & \centering \multirow{2}{*}{Twomey} \tabularnewline
& & \centering fin. diff. & \tabularnewline
\hline
\hline
constrained & \centering 58.7991 & \centering 58.7264 & \centering 58.7173 \tabularnewline
\hline
Morozov & \centering 525.3771 & \centering 557.9417 & \centering 579.4758 \tabularnewline
\hline
unconstrained & \centering 70.6745 & \centering 64.0843 & \centering 67.1404 \tabularnewline
\hline
BIC & \centering $1.0245 \cdot 10^{4}$ & \centering $1.0246 \cdot 10^{4}$ & \centering $1.0248 \cdot 10^{4}$ \tabularnewline
\hline
\end{tabular}
\end{table}

\newpage

\begin{table}[h!]
\centering
\begin{tabular}{||m{2.5cm}||m{2.5cm}|m{2.5cm}|m{2.5cm}||}
\hline
\multicolumn{4}{||c||}{\multirow{3}{*}{}} \tabularnewline[-2.5ex]
\multicolumn{4}{||c||}{RRSB Distribution} \tabularnewline
\multicolumn{4}{||c||}{} \tabularnewline[-2.5ex]
\hline
\multirow{4}{*}{\vspace{0.2cm}method} & \multicolumn{3}{c||}{} \tabularnewline[-2.5ex]
 & \multicolumn{3}{c||}{worst case $L^{2}$-errors (\%)} \tabularnewline
& \multicolumn{3}{c||}{} \tabularnewline[-2.5ex]
\cline{2-4}
& \centering \multirow{2}{*}{Tikhonov} & \centering min. first & \centering \multirow{2}{*}{Twomey} \tabularnewline
& & \centering fin. diff. & \tabularnewline
\hline
\hline
constrained & \centering 118.8957  & \centering 114.6170 & \centering 114.6434 \tabularnewline
\hline
Morozov & \centering 389.9993 & \centering 379.7349 & \centering 399.7682 \tabularnewline
\hline
unconstrained & \centering 85.2596 & \centering 81.1347 & \centering 82.4548 \tabularnewline
\hline
BIC & \centering $2.2651 \cdot 10^{4}$ & \centering $2.2656 \cdot 10^{4}$ & \centering $2.2663 \cdot 10^{4}$ \tabularnewline
\hline
\end{tabular}
\end{table}

\begin{table}[h!]
\centering
\begin{tabular}{||m{2.5cm}||m{2.5cm}|m{2.5cm}|m{2.5cm}||}
\hline
\multicolumn{4}{||c||}{\multirow{3}{*}{}} \tabularnewline[-2.5ex]
\multicolumn{4}{||c||}{Hedrih Distribution} \tabularnewline
\multicolumn{4}{||c||}{} \tabularnewline[-2.5ex]
\hline
\multirow{4}{*}{\vspace{0.2cm}method} & \multicolumn{3}{c||}{} \tabularnewline[-2.5ex]
 & \multicolumn{3}{c||}{worst case $L^{2}$-errors (\%)} \tabularnewline
& \multicolumn{3}{c||}{} \tabularnewline[-2.5ex]
\cline{2-4}
& \centering \multirow{2}{*}{Tikhonov} & \centering min. first & \centering \multirow{2}{*}{Twomey} \tabularnewline
& & \centering fin. diff. & \tabularnewline
\hline
\hline
constrained & \centering 51.1577 & \centering 38.2242 & \centering 35.3611 \tabularnewline
\hline
Morozov & \centering 200.3302 & \centering 224.3174 & \centering 239.6155 \tabularnewline
\hline
unconstrained & \centering 56.6067 & \centering 56.5325 & \centering 56.4147 \tabularnewline
\hline
BIC & \centering $1.0273 \cdot 10^{4}$ & \centering $1.0273 \cdot 10^{4}$ & \centering $1.0273 \cdot 10^{4}$ \tabularnewline
\hline
\end{tabular}
\end{table}

\subsubsection{Worst Case Run Times}

\begin{table}[h!]
\centering
\begin{tabular}{||m{2.5cm}||m{2.5cm}|m{2.5cm}|m{2.5cm}||}
\hline
\multicolumn{4}{||c||}{\multirow{3}{*}{}} \tabularnewline[-2.5ex]
\multicolumn{4}{||c||}{Log-Normal Distribution} \tabularnewline
\multicolumn{4}{||c||}{} \tabularnewline[-2.5ex]
\hline
\multirow{4}{*}{\vspace{0.2cm}method} & \multicolumn{3}{c||}{} \tabularnewline[-2.5ex]
 & \multicolumn{3}{c||}{worst case run times (s)} \tabularnewline
& \multicolumn{3}{c||}{} \tabularnewline[-2.5ex]
\cline{2-4}
& \centering \multirow{2}{*}{Tikhonov} & \centering min. first & \centering \multirow{2}{*}{Twomey} \tabularnewline
& & \centering fin. diff. & \tabularnewline
\hline
\hline
constrained & \centering 6.2403 & \centering 7.5913 & \centering 7.0653 \tabularnewline
\hline
Morozov & \centering 2.9394 & \centering 3.3588 & \centering 4.3115 \tabularnewline
\hline
unconstrained & \centering 0.7136 & \centering 0.7075 & \centering 0.6234 \tabularnewline
\hline
BIC & \centering 0.1126 & \centering 0.1689 & \centering 0.1645 \tabularnewline
\hline
\end{tabular}
\end{table}

\begin{table}[h!]
\centering
\begin{tabular}{||m{2.5cm}||m{2.5cm}|m{2.5cm}|m{2.5cm}||}
\hline
\multicolumn{4}{||c||}{\multirow{3}{*}{}} \tabularnewline[-2.5ex]
\multicolumn{4}{||c||}{RRSB Distribution} \tabularnewline
\multicolumn{4}{||c||}{} \tabularnewline[-2.5ex]
\hline
\multirow{4}{*}{\vspace{0.2cm}method} & \multicolumn{3}{c||}{} \tabularnewline[-2.5ex]
 & \multicolumn{3}{c||}{worst case run times (s)} \tabularnewline
& \multicolumn{3}{c||}{} \tabularnewline[-2.5ex]
\cline{2-4}
& \centering \multirow{2}{*}{Tikhonov} & \centering min. first & \centering \multirow{2}{*}{Twomey} \tabularnewline
& & \centering fin. diff. & \tabularnewline
\hline
\hline
constrained & \centering 9.6938 & \centering 10.1606 & \centering 10.0524 \tabularnewline
\hline
Morozov & \centering 2.2293 & \centering 2.6416 & \centering 2.3914 \tabularnewline
\hline
unconstrained & \centering 0.4096 & \centering 0.4453 & \centering 0.4661 \tabularnewline
\hline
BIC & \centering 0.0925 & \centering 0.0896 & \centering 0.1010 \tabularnewline
\hline
\end{tabular}
\end{table}

\newpage

\begin{table}[h!]
\centering
\begin{tabular}{||m{2.5cm}||m{2.5cm}|m{2.5cm}|m{2.5cm}||}
\hline
\multicolumn{4}{||c||}{\multirow{3}{*}{}} \tabularnewline[-2.5ex]
\multicolumn{4}{||c||}{Hedrih Distribution} \tabularnewline
\multicolumn{4}{||c||}{} \tabularnewline[-2.5ex]
\hline
\multirow{4}{*}{\vspace{0.2cm}method} & \multicolumn{3}{c||}{} \tabularnewline[-2.5ex]
 & \multicolumn{3}{c||}{worst case run times (s)} \tabularnewline
& \multicolumn{3}{c||}{} \tabularnewline[-2.5ex]
\cline{2-4}
& \centering \multirow{2}{*}{Tikhonov} & \centering min. first & \centering \multirow{2}{*}{Twomey} \tabularnewline
& & \centering fin. diff. & \tabularnewline
\hline
\hline
constrained & \centering 2.9746 & \centering 2.7281 & \centering 2.7897 \tabularnewline
\hline
Morozov & \centering 1.7559 & \centering 1.8946 & \centering 2.0438 \tabularnewline
\hline
unconstrained & \centering 0.3564 & \centering 0.4351 & \centering 0.3860 \tabularnewline
\hline
BIC & \centering 0.0556 & \centering 0.0625 & \centering 0.0836 \tabularnewline
\hline
\end{tabular}
\end{table}

\subsection{Conclusion}

The constrained method had the smallest average $L^{2}$-errors and close to zero failure rates. Only for the RRSB distributions were two, one, and one failures out of $1000$ inversions recorded for the different priors, respectively. The overall worst case reconstruction error of $118.8957$\% was only moderately above $100$\%. It needed the longest run times from all methods, but even the overall worst case run time of $10.9830$ seconds was clearly below our thirty-second requirement. The difference of the average $L^{2}$-errors depending on the three priors we applied was not very prominent. For RRSB and Hedrih distributions they seem to decrease by ca. $0.5$ to $1$\% from the Tikhonov to the minimal first finite differences to the Twomey priors, whereas for log-normal distributions the opposite behavior is the case. For the other inversion methods the $L^{2}$-errors behave similarly depending on the priors. Therefore we cannot determine a prior out of the three we used which always yields the smallest average $L^{2}$-error. 

For the Morozov method the average $L^{2}$-errors were for log-normal and RRSB distributions about $4$ to $6$\% higher compared to those of the constrained method, but for Hedrih distributions the errors were $11$\% larger. The average run times represented only about one fifth of those of the constrained method.  However, the numbers of failures was significantly higher. For log-normal and Hedrih distributions roughly $3$\% of all inversions failed, but for RRSB distributions this was up to $6$\%. The overall worst case $L^{2}$-error of $579.4758$\% was clearly higher than $100$\%.

The run times of the unconstrained method were one third to one half of the Morozov method run times. The unconstrained method was the only method without any failures. The overall worst case $L^{2}$-error was a relatively moderate $85.2596$\%, but the average $L^{2}$-errors were $5$ to $12$\% bigger than the Morozov method $L^{2}$-errors and already $1.5$ to $3$ times as big as the constrained method $L^{2}$-errors.

The BIC method was by far the fastest one with run times of only a few hundredths of a second, but the average $L^{2}$ errors ranging from ca. $40$ to $80$\% were rather poor. The overall worst case $L^{2}$-error was even $2.2663 \cdot 10^{4}$\%.
 
For practical FASP experiments we conclude that the constrained method performed best, because its average $L^{2}$-errors were smallest, had virtually no failures, and clearly satisfied our thirty-seconds run-time limit even in the worst cases.

\section{Two-Component Aerosols}

In the preceding sections it was assumed that the aerosol particles consist of a known material, and therefore the refractive indices $m_{part}(l)$ needed to compute the extinction efficiency $Q_{ext}(m_{med}(l), m_{part}(l), r, l)$ were given exactly as well. But this is not generally the case in real experimental measurements where typically both size distributions and optical properties of scattering particles are unknown. In the ideal case we could set up an additional device for measuring the aerosol refractive indices and perform a two-stage measurement process, where the first step is to retrieve the refractive indices as preparation for the second step of reconstructing the size distribution, but this is not practical. Indeed all measurement techniques for optical properties of aerosol particles demand a pretreatment of the aerosol itself such as vaporizing it into its gas phase or transforming it into a monodisperse aerosol. This would make the FASP too inefficient to be of practical use.

In real applications we simply want to examine some aerosol components of particular interest. Thus we assume that the aerosol to be investigated is a mixture of a small number of known materials, such that only the problem remains to retrieve the volume fractions of these materials in the whole composite aerosol. As an initial explorative step into this general problem we further assume that the aerosol is made up of only two materials.

To compute the refractive indices of composite aerosols from those of their pure components so-called \textit{mixing rules} are used.  Some of these are compared in \cite{SPPV06}. Let $m_1 = k_1 + i n_1$ and $m_2 = k_2 + i n_2$ be the refractive indices of two aerosol components for a wavelength $l$ of the incident light. We adopt the most commonly used rule, the \textit{Lorentz-Lorenz rule}. Here the total refractive $m_{tot} = k_{tot} + i n_{tot}$ is obtained from the relation
\begin{equation}
\frac{m_{tot}^2 - 1}{m_{tot}^2 + 2} = f_1 \frac{m_{1}^2 - 1}{m_{1}^2 + 2} + f_2 \frac{m_{2}^2 - 1}{m_{2}^2 + 2},
\label{Lorentz-Lorenz}
\end{equation} 
where $f_1$ and $f_2$ are the volume fractions of the components.

Now our new problem is to invert the parameter-dependent integral equation
\begin{equation}
\int_{r_{min}}^{r_{max}} k_{p}(r, l)n(r) dr = e(l),
\end{equation}
where the sought-after parameter $p \in [0, 1]$ characterizes the unknown volume fractions. Let $m_{p}(l)$ denote the solution $m_{tot}$ of \eqref{Lorentz-Lorenz} with $f_1 = p$ and $f_2 = 1 - p$. Then the $p$-dependent kernel function is given by
\begin{equation*}
k_{p}(r, l) = \pi r^{2} Q_{ext}(m_{med}(l), m_{p}(l), r, l).
\end{equation*}
Mathematically this means that in addition to inverting it we have to identify the ``right'' integral operator $K_{p}$ from the set 
\begin{equation*}
\left\{(K_{p}n)(l) := \int_{r_{min}}^{r_{max}} k_{p}(r,l) n(r) dr \; \bigg| \; p \in [0, 1] \right\}.
\end{equation*}

We can easily check that $k_{p}(r, l)$ depends continuously on $p$ and therefore so do the discrete approximations $\boldsymbol{K}_{k, p}$ to $K_{p}$ as well. We again make Assumption \ref{covariance}. So by setting 
\begin{equation*}
\boldsymbol{K}_{p} := \boldsymbol{\Sigma}^{-\frac{1}{2}}\boldsymbol{K}_{k, p} \quad \text{and} \quad \boldsymbol{r} = \boldsymbol{\Sigma}^{-\frac{1}{2}}(\boldsymbol{e}_{true} + \boldsymbol{\delta})
\end{equation*}
as in Section \ref{convergence1} we obtain the $p$-parametrized quadratic programming problem
\begin{equation}
\min_{\boldsymbol{n} \in \mathbb{R}^{N}} \textstyle{\frac{1}{2}}\|\boldsymbol{K}_{p}\boldsymbol{n} - \boldsymbol{r}\|_{2}^{2} + \textstyle{\frac{1}{2}}\gamma\|\boldsymbol{n}\|_{2}^{2} \quad \text{s.t.} \quad \boldsymbol{Cn} \leq \boldsymbol{b}.
\label{quadprog_p}
\end{equation}
as in Section \ref{motivation} for the computation of the maximum a posteriori solution.

\subsection{Fraction Retrieval for two Aerosol Components}

For the determination of the parameter $p$ we modify the adaptive model-generation algorithm from Section \ref{section_model_generation}. As a preparation we prove a continuity result.
\begin{prp}
\textit{The minimizer $\boldsymbol{n}_{p}$ of \eqref{quadprog_p} for $\gamma = 0$ depends continuously on the kernel matrix $\boldsymbol{K}_{p}$.}
\end{prp}

\begin{prf}
Let $p_1, p_2 \in [0, 1]$ be arbitrary. We write 
\begin{equation*}
\boldsymbol{K}_{p_1} =: \boldsymbol{K} \quad \text{and} \quad \boldsymbol{K}_{p_2} =: \boldsymbol{K} + \boldsymbol{S},
\end{equation*}
hence $\boldsymbol{S} = \boldsymbol{K}_{p_2} - \boldsymbol{K}_{p_1}$. From the continuous dependence of $\boldsymbol{K}_{p}$ on $p$ we have 
\begin{equation}
\lim_{p_2 \to p_1} \boldsymbol{S} = 0.
\label{S_lim}
\end{equation}

The first-order necessary conditions for the minimizers $\boldsymbol{n}_{p_1}$ and $\boldsymbol{n}_{p_2}$ of \eqref{quadprog_p} for $p = p_1$ and $p = p_2$ are given by the relations
\begin{align}
& \boldsymbol{K}^{T}\boldsymbol{K}\boldsymbol{n}_{p_1} - \boldsymbol{K}^{T}\boldsymbol{r} + \boldsymbol{C}^{T}\boldsymbol{q}_{p_1} = 0 \label{KKT1} \\
\text{and} \quad 
(&\boldsymbol{K}^{T}\boldsymbol{K}  + \boldsymbol{K}^{T}\boldsymbol{S} +  \boldsymbol{S}^{T}\boldsymbol{K} +  \boldsymbol{S}^{T}\boldsymbol{S})\boldsymbol{n}_{p_2} - (\boldsymbol{K} +\boldsymbol{S})^{T}\boldsymbol{r} + \boldsymbol{C}^{T}\boldsymbol{q}_{p_2} = 0,
\label{KKT2}
\end{align}
with vectors $\boldsymbol{q}_{p_1} \geq 0$, $\boldsymbol{q}_{p_2} \geq 0$. Subtracting \eqref{KKT1} from \eqref{KKT2} yields
\begin{equation*}
\boldsymbol{K}^{T}\boldsymbol{K}(\boldsymbol{n}_{p_2} - \boldsymbol{n}_{p_1})  + (\boldsymbol{K}^{T}\boldsymbol{S} +  \boldsymbol{S}^{T}\boldsymbol{K} +  \boldsymbol{S}^{T}\boldsymbol{S})\boldsymbol{n}_{p_2} - \boldsymbol{S}^{T}\boldsymbol{r} + \boldsymbol{C}^{T}(\boldsymbol{q}_{p_2} - \boldsymbol{q}_{p_1}) = 0.
\end{equation*}
Forming the scalar product with $\boldsymbol{n}_{p_2} - \boldsymbol{n}_{p_1}$ yields
\begin{equation*}
\begin{split}
& \langle \boldsymbol{n}_{p_2} - \boldsymbol{n}_{p_1}, \boldsymbol{K}^{T}\boldsymbol{K}(\boldsymbol{n}_{p_2} - \boldsymbol{n}_{p_1})\rangle + \langle \boldsymbol{n}_{p_2} - \boldsymbol{n}_{p_1}, (\boldsymbol{K}^{T}\boldsymbol{S} +  \boldsymbol{S}^{T}\boldsymbol{K} +  \boldsymbol{S}^{T}\boldsymbol{S})\boldsymbol{n}_{p_2}\rangle \\
- \; & \langle\boldsymbol{n}_{p_2} - \boldsymbol{n}_{p_1}, \boldsymbol{S}^{T}\boldsymbol{r}\rangle + \langle\boldsymbol{C}(\boldsymbol{n}_{p_2} - \boldsymbol{n}_{p_1}), \boldsymbol{q}_{p_2} - \boldsymbol{q}_{p_1}\rangle = 0.
\end{split}
\end{equation*}
With \eqref{S_lim} we obtain in the limit $p_2 \rightarrow p_1$
\begin{equation*}
\lim_{p_2 \to p_1}\Big( \langle \boldsymbol{n}_{p_2} - \boldsymbol{n}_{p_1}, \boldsymbol{K}^{T}\boldsymbol{K}(\boldsymbol{n}_{p_2} - \boldsymbol{n}_{p_1})\rangle + \langle\boldsymbol{C}(\boldsymbol{n}_{p_2} - \boldsymbol{n}_{p_1}), \boldsymbol{q}_{p_2} - \boldsymbol{q}_{p_1}\rangle \Big) = 0.
\end{equation*}
A calculation as in the proof of Lemma \ref{monotonicity} shows
\begin{equation*}
\langle\boldsymbol{C}(\boldsymbol{n}_{p_2} - \boldsymbol{n}_{p_1}), \boldsymbol{q}_{p_2} - \boldsymbol{q}_{p_1}\rangle \geq 0,
\end{equation*}
which finally implies
\begin{equation*}
\lim_{p_2 \to p_1}\boldsymbol{n}_{p_2} = \boldsymbol{n}_{p_1}.
\end{equation*}
\hfill $\square$
\end{prf}

From the last proposition we directly obtain an existence result for an optimal $p$.
\begin{crl}
\textit{For $\gamma = 0$ the residual $\|\boldsymbol{K}_{p}\boldsymbol{n}_{p} - \boldsymbol{r}\|_2$ of the minimizer of \eqref{quadprog_p} depends continuously on $p$, so there exists a $p \in [0, 1]$ for which it attains its minimal value.}
\hfill $\square$
\end{crl}

Our next step is to find a condition for uniqueness of this minimizer for $\gamma = 0$.

\begin{prp}
\label{unique_fraction}
\textit{Let $\gamma = 0$ and $p \in [0, 1]$ be such that $\|\boldsymbol{K}_{p}\boldsymbol{n}_{p} - \boldsymbol{r}\|_2^{2}$ minimizes all Tikhonov functionals in \eqref{quadprog_p} over the parameter range $[0, 1]$. Lets $s \in [0, 1]$, $s \neq p$, be arbitrary. Then if}
\begin{equation}
\langle \boldsymbol{K}_{p}\boldsymbol{n}_{p} - \boldsymbol{K}_{s}\boldsymbol{n}_{s}, \boldsymbol{r} \rangle \neq 0
\label{cond_unique}
\end{equation}
\textit{holds, the minimizing parameter $p$ is unique.}
\end{prp} 

\begin{prf}
The necessary conditions for $\boldsymbol{n}_{p}$ and $\boldsymbol{n}_{s}$ to be a minimizer of \eqref{quadprog_p} are given by
\begin{align}
\boldsymbol{K}_{p}^{T}\boldsymbol{K}_{p}\boldsymbol{n}_{p} - \boldsymbol{K}_{p}^{T}\boldsymbol{r} + \boldsymbol{C}^{T}\boldsymbol{q}_{p} & = 0 \label{KKT_unique_1} \\
\text{and} \quad \boldsymbol{K}_{s}^{T}\boldsymbol{K}_{s}\boldsymbol{n}_{s} - \boldsymbol{K}_{s}^{T}\boldsymbol{r} + \boldsymbol{C}^{T}\boldsymbol{q}_{s} & = 0
\label{KKT_unique_2}
\end{align}
with vectors $\boldsymbol{q}_{p} \geq 0$, $\boldsymbol{q}_{s} \geq 0$. Assume
\begin{equation*}
\|\boldsymbol{K}_{p}\boldsymbol{n}_{p} - \boldsymbol{r}\|_2^{2} = \|\boldsymbol{K}_{s}\boldsymbol{n}_{s} - \boldsymbol{r}\|_2^{2},
\end{equation*}
which is equivalent to
\begin{equation}
\langle \boldsymbol{n}_{p}, \boldsymbol{K}_{p}^{T}\boldsymbol{K}_{p}\boldsymbol{n}_{p}\rangle - \langle \boldsymbol{n}_{s}, \boldsymbol{K}_{s}^{T}\boldsymbol{K}_{s}\boldsymbol{n}_{s}\rangle = 2\langle \boldsymbol{K}_{p}\boldsymbol{n}_{p} - \boldsymbol{K}_{s}\boldsymbol{n}_{s}, \boldsymbol{r} \rangle. 
\label{intermediate_step}
\end{equation}

We form the scalar products of \eqref{KKT_unique_1} with $\boldsymbol{n}_{p}$ and of \eqref{KKT_unique_2} with $\boldsymbol{n}_{s}$. Then forming the difference of the resulting equations gives
\begin{equation*}
\langle \boldsymbol{n}_{p}, \boldsymbol{K}_{p}^{T}\boldsymbol{K}_{p}\boldsymbol{n}_{p}\rangle - \langle \boldsymbol{n}_{s}, \boldsymbol{K}_{s}^{T}\boldsymbol{K}_{s}\boldsymbol{n}_{s}\rangle - \langle \boldsymbol{K}_{p}\boldsymbol{n}_{p} - \boldsymbol{K}_{s}\boldsymbol{n}_{s}, \boldsymbol{r} \rangle + \langle\boldsymbol{C}\boldsymbol{n}_{p}, \boldsymbol{q}_{p}\rangle - \langle \boldsymbol{C}\boldsymbol{n}_{s}, \boldsymbol{q}_{s}\rangle = 0.
\end{equation*}
Inserting \eqref{intermediate_step} yields
\begin{equation*}
\langle \boldsymbol{K}_{p}\boldsymbol{n}_{p} - \boldsymbol{K}_{s}\boldsymbol{n}_{s}, \boldsymbol{r} \rangle + \langle \boldsymbol{C}\boldsymbol{n}_{p}, \boldsymbol{q}_{p}\rangle - \langle \boldsymbol{C}\boldsymbol{n}_{s}, \boldsymbol{q}_{s}\rangle = 0.
\end{equation*}
From \eqref{KKT_lin_constr_4} with $\boldsymbol{b} = 0$ we conclude $\langle \boldsymbol{Cn}_{p}, \boldsymbol{q}_{p} \rangle = 0$ and analogously $\langle \boldsymbol{C}\boldsymbol{n}_{s}, \boldsymbol{q}_{s}\rangle = 0$. But then
\begin{equation*}
\langle \boldsymbol{K}_{p}\boldsymbol{n}_{p} - \boldsymbol{K}_{s}\boldsymbol{n}_{s}, \boldsymbol{r} \rangle = 0,
\end{equation*}
which contradicts \eqref{cond_unique}. Thus if \eqref{cond_unique} holds, the minimizing parameter $p$ must be unique. 

\hfill $\square$
\end{prf}

Condition \eqref{cond_unique} demands that the kernel matrices $\boldsymbol{K}_{s}$ are sufficiently different to $\boldsymbol{K}_{p}$ so that we get distinguishable residuals of the unregularized solutions. If an $s \in [0, 1]$ happens to exist with $\boldsymbol{K}_{s} = \boldsymbol{K}_{p}$, condition \eqref{cond_unique} cannot be fulfilled.
Unfortunately we are not currently able to check this condition a priori.

We conclude this section with investigating how the unregularized residuals behave for moderate noise levels. In the following the superscript $\delta$ marks solutions of \eqref{quadprog_p} for a data vector $\boldsymbol{r}$ contaminated with noise.

\begin{prp}
\label{noise_limit}
\textit{We assume that condition \eqref{cond_unique} holds for any noise vector satisfying $0 \leq \|\boldsymbol{\delta}\|_{2} \leq \delta$. Let $\boldsymbol{n}_{t}$ be the minimizer for the true noise-free model, i.e.\ the parameter $t \in [0, 1]$ yields the minimal residual $\|\boldsymbol{K}_{t}\boldsymbol{n}_{t} - \boldsymbol{\Sigma}^{-\frac{1}{2}}\boldsymbol{e}_{true}\|_{2}^{2}$ over the whole parameter interval $[0, 1]$. Let $p = p(\boldsymbol{\delta}) \in [0, 1]$ be the parameter yielding the minimal unregularized residual $\|\boldsymbol{K}_{p}\boldsymbol{n}_{p}^{\delta} - \boldsymbol{\Sigma}^{-\frac{1}{2}}(\boldsymbol{e}_{true} + \boldsymbol{\delta})\|_{2}^{2}$ for the noisy data vector $\boldsymbol{\Sigma}^{-\frac{1}{2}}(\boldsymbol{e}_{true} + \boldsymbol{\delta})$. Then holds $\lim_{\|\boldsymbol{\delta}\|_{2} \to 0} p(\boldsymbol{\delta}) = t$.}
\end{prp}

\begin{prf}
To shorten notation we write $\boldsymbol{r}_{true} := \boldsymbol{\Sigma}^{-\frac{1}{2}}\boldsymbol{e}_{true}$ and $\boldsymbol{\rho} := \boldsymbol{\Sigma}^{-\frac{1}{2}}\boldsymbol{\delta}$. Let $\boldsymbol{n}_{t}^{\delta}$ be the minimizer for the parameter $t$ and the noisy data vector $\boldsymbol{\Sigma}^{-\frac{1}{2}}(\boldsymbol{e}_{true} + \boldsymbol{\delta})$, i.e.\ 
\begin{equation*}
\boldsymbol{n}_{t}^{\delta} = \underset{\boldsymbol{n} \in \mathbb{R}^{N}}{\mathrm{argmin}} \; \textstyle{\frac{1}{2}}\|\boldsymbol{K}_{t}\boldsymbol{n} - (\boldsymbol{r}_{true} + \boldsymbol{\rho})\|_2^{2} \quad \text{s.t.} \quad \boldsymbol{Cn} \leq \boldsymbol{b}.
\end{equation*}
Then we have the first order necessary conditions
\begin{align}
\boldsymbol{K}_{t}^{T}\boldsymbol{K}_{t}\boldsymbol{n}_{t} - \boldsymbol{K}_{t}^{T}\boldsymbol{r}_{true} + \boldsymbol{C}\boldsymbol{q}_{t} & = 0 \label{noise-free} \\
\boldsymbol{K}_{t}^{T}\boldsymbol{K}_{t}\boldsymbol{n}_{t}^{\delta} - \boldsymbol{K}_{t}^{T}(\boldsymbol{r}_{true} + \boldsymbol{\rho}) + \boldsymbol{C}\boldsymbol{q}_{t}^{\delta} & = 0, \label{noisy}
\end{align}
with vectors $\boldsymbol{q}_{t}^{\delta} \geq 0$, $\boldsymbol{q}_{t} \geq 0$.
Subtracting \eqref{noisy} from \eqref{noise-free} and scalar multiplying the result with $\boldsymbol{n}_{t} - \boldsymbol{n}_{t}^{\delta}$ gives
\begin{equation*}
\langle \boldsymbol{n}_{t} - \boldsymbol{n}_{t}^{\delta}, \boldsymbol{K}_{t}^{T}\boldsymbol{K}_{t}(\boldsymbol{n}_{t} - \boldsymbol{n}_{t}^{\delta})\rangle + \langle \boldsymbol{n}_{t} - \boldsymbol{n}_{t}^{\delta}, \boldsymbol{K}_{t}^{T}\boldsymbol{\rho}\rangle + \langle \boldsymbol{C}(\boldsymbol{n}_{t} - \boldsymbol{n}_{t}^{\delta}), \boldsymbol{q}_{t} - \boldsymbol{q}_{t}^{\delta}\rangle = 0.
\end{equation*}
As in the proof of Lemma \ref{monotonicity} this establishes 
\begin{align*}
\|\boldsymbol{K}_{t}(\boldsymbol{n}_{t} - \boldsymbol{n}_{t}^{\delta})\|_{2}^{2} & \leq \langle \boldsymbol{n}_{t}^{\delta} - \boldsymbol{n}_{t},\boldsymbol{K}_{t}^{T}\boldsymbol{\rho}\rangle \\
& \leq \|\boldsymbol{K}_{t}(\boldsymbol{n}_{t} - \boldsymbol{n}_{t}^{\delta})\|_{2}\|\boldsymbol{\Sigma}^{-\frac{1}{2}}\boldsymbol{\delta}\|_{2},
\end{align*}
which gives
\begin{equation*}
\|\boldsymbol{K}_{t}(\boldsymbol{n}_{t} - \boldsymbol{n}_{t}^{\delta})\|_{2} = \mathcal{O}(\|\boldsymbol{\delta}\|_{2}).
\end{equation*}

Now since the parameter $p$ minimizes the residuals for the noisy vector $\boldsymbol{r}_{true} + \boldsymbol{\rho}$ we can estimate
\begin{align*}
\|\boldsymbol{K}_{p}\boldsymbol{n}_{p}^{\delta} - (\boldsymbol{r}_{true} + \boldsymbol{\rho})\|_{2} & \leq \|\boldsymbol{K}_{t}\boldsymbol{n}_{t}^{\delta} - (\boldsymbol{r}_{true} + \boldsymbol{\rho})\|_{2} \\
& \leq \|\boldsymbol{K}_{t}\boldsymbol{n}_{t} - \boldsymbol{r}_{true}\|_{2} + \|\boldsymbol{K}_{t}(\boldsymbol{n}_{t}^{\delta} - \boldsymbol{n}_{t})\|_{2} + \|\boldsymbol{\Sigma}^{-\frac{1}{2}}\boldsymbol{\delta}\|_{2}.
\end{align*}
With the previous finding we see that the upper bound in the last inequality converges to the residual $\|\boldsymbol{K}_{t}\boldsymbol{n}_{t} - \boldsymbol{r}_{true}\|_{2}$ for $\|\boldsymbol{\delta}\|_{2} \rightarrow 0$. Thus we obtain in the limit
\begin{equation*}
\lim_{\|\boldsymbol{\delta}\|_{2} \to 0}\|\boldsymbol{K}_{p}\boldsymbol{n}_{p}^{\delta} - (\boldsymbol{r}_{true} + \boldsymbol{\rho})\|_{2} \leq \|\boldsymbol{K}_{t}\boldsymbol{n}_{t} - \boldsymbol{r}_{true}\|_{2}.
\end{equation*}
By definition of $t$ we have
\begin{equation*}
\|\boldsymbol{K}_{t}\boldsymbol{n}_{t} - \boldsymbol{r}_{true}\|_{2} \leq \lim_{\|\boldsymbol{\delta}\|_{2} \to 0}\|\boldsymbol{K}_{p}\boldsymbol{n}_{p}^{\delta} - (\boldsymbol{r}_{true} + \boldsymbol{\rho})\|_{2},
\end{equation*}
so condition \eqref{cond_unique} finally implies $\lim_{\|\boldsymbol{\delta}\|_{2} \to 0} p(\boldsymbol{\delta}) = t$.
\hfill $\square$
\end{prf}

\subsection{Convergence Analysis}

In this section we show that the regularized solutions from the retrieved aerosol fraction converge to the true solution from the true fraction as the noise level approaches zero. This means that we generalize Theorem \ref{conv_thm1} to the case where the underlying true linear operator must be identified from a known set of possible operators. 

\begin{thm}
\textit{Under Assumption \ref{covariance}, if condition \ref{cond_unique} is satisfied for all noise vectors $\boldsymbol{\delta}$ of random variables, then we have for any $\alpha(\delta)$ with the properties $\lim_{\delta \to 0} \: \alpha(\delta) = 0$ and $\lim_{\delta \to 0} \: \frac{\delta^{2}}{\alpha(\delta)} = 0$ that $\lim_{\delta \to 0} \: \mathbb{E}\big(\|\boldsymbol{n}_{p}^{\delta, \alpha(\delta)} - \boldsymbol{n}_{t}\|_{2}\big) = 0$.  Here $\boldsymbol{n}_{p}^{\delta, \alpha(\delta)}$ is the regularized solution for the retrieved fraction parameter $p = p(\boldsymbol{\delta})$.}
\end{thm}

\begin{prf}
We again use the notations $\boldsymbol{r}_{true} := \boldsymbol{\Sigma}^{-\frac{1}{2}}\boldsymbol{e}_{true}$ and $\boldsymbol{\rho} := \boldsymbol{\Sigma}^{-\frac{1}{2}}\boldsymbol{\delta}$. Let $p = p(\boldsymbol{\delta})$ be the fraction parameter retrieved by minimizing the unregularized residuals. We write

\begin{align*}
\boldsymbol{n}_{p}^{\delta, \alpha} & := \underset{\boldsymbol{n} \: \in \: \mathbb{R}^{N}}{\mathrm{argmin}} \; \textstyle{\frac{1}{2}}\|\boldsymbol{K}_{p}\boldsymbol{n} - (\boldsymbol{r}_{true} + \boldsymbol{\rho})\|_{2}^{2} + \textstyle{\frac{1}{2}}\alpha\|\boldsymbol{n}\|_{2}^{2} \quad \text{s.t.} \quad \boldsymbol{Cn} \leq \boldsymbol{b},\\
\boldsymbol{n}_{p}^{\alpha} & := \underset{\boldsymbol{n} \: \in \: \mathbb{R}^{N}}{\mathrm{argmin}} \; \textstyle{\frac{1}{2}}\|\boldsymbol{K}_{p}\boldsymbol{n} - \boldsymbol{r}_{true}\|_{2}^{2} + \textstyle{\frac{1}{2}}\alpha\|\boldsymbol{n}\|_{2}^{2} \quad \text{s.t.} \quad \boldsymbol{Cn} \leq \boldsymbol{b},\\
\boldsymbol{n}_{p} & := \underset{\boldsymbol{n} \: \in \: \mathbb{R}^{N}}{\mathrm{argmin}} \; \textstyle{\frac{1}{2}}\|\boldsymbol{K}_{p}\boldsymbol{n} - \boldsymbol{r}_{true}\|_{2}^{2} \quad \text{s.t.} \quad \boldsymbol{Cn} \leq \boldsymbol{b},\\
\boldsymbol{n}_{t} & := \underset{\boldsymbol{n} \: \in \: \mathbb{R}^{N}}{\mathrm{argmin}} \; \textstyle{\frac{1}{2}}\|\boldsymbol{K}_{t}\boldsymbol{n} - \boldsymbol{r}_{true}\|_{2}^{2} \quad \text{s.t.} \quad \boldsymbol{Cn} \leq \boldsymbol{b}.
\end{align*}
Then we have the estimate
\begin{equation*}
\mathbb{E}\big(\|\boldsymbol{n}_{p}^{\delta, \alpha} - \boldsymbol{n}_{t}\|_{2}\big) \leq \mathbb{E}\big(\|\boldsymbol{n}_{p}^{\delta, \alpha} - \boldsymbol{n}_{p}^{\alpha}\|_{2}\big) + \mathbb{E}\big(\|\boldsymbol{n}_{p}^{\alpha} - \boldsymbol{n}_{p}\|_{2}\big) + \mathbb{E}(\|\boldsymbol{n}_{p} - \boldsymbol{n}_{t}\|_{2}).
\end{equation*}
For the first term in the upper bound, the estimate
\begin{equation*}
\mathbb{E}\big(\|\boldsymbol{n}_{p}^{\delta, \alpha} - \boldsymbol{n}_{p}^{\alpha}\|_{2}\big) = \frac{\mathcal{O}(\delta)}{\alpha^{\frac{1}{2}}}
\end{equation*}
follows from Proposition \ref{bounds}. For the second term, Proposition \ref{convergence} gives
\begin{equation*}
\lim_{\alpha \to 0} \; \mathbb{E}\big(\|\boldsymbol{n}_{p}^{\alpha} - \boldsymbol{n}_{p}\|_{2}\big) = 0.
\end{equation*}
Finally, for the third term from Proposition \ref{noise_limit} follows
\begin{equation*}
\lim_{\delta \to 0} \; \mathbb{E}(\|\boldsymbol{n}_{p} - \boldsymbol{n}_{t}\|_{2}) = 0.
\end{equation*}
This altogether proves our claim.
\hfill $\square$
\end{prf}

\subsection{Model Generation under Nonnegativity Constraints for Two-Component Aerosols}
\label{model_gen_2comp}

Proposition \ref{unique_fraction} motivates us to use the unregularized residuals as model generation criterion, which means that we determine those parameters $s$, where they are small. In presence of moderate measurement noise in $\boldsymbol{e}$ these parameters lie in the vicinity of the unique true parameter $p$ as was shown in the proof of Proposition \ref{noise_limit}. In the following we discuss the model generation algorithm extended for two-component aerosols. As in Section \ref{section_model_generation} we compute collocation grids with $N_{1} < ... < N_{m}$ points and select a grid of Morozov safety factors $\tau_{1} < ... < \tau_{s}$. Furthermore the refractive indices $k_{1}(l_{1}) + i n_{1}(l_{1}), ..., k_{1}(l_{N_{l}}) + i n_{1}(l_{N_{l}})$ and $k_{2}(l_{1}) + i n_{2}(l_{1}), ..., k_{2}(l_{N_{l}}) + i n_{2}(l_{N_{l}})$ of two pure aerosol components depending on wavelengths $l_{1}, ..., l_{N_{l}}$ are given.

\newpage

\begin{algorithm}[H]
\caption{Model Generation for Two-Component Aerosols}
\label{model_generation_2comp}
\begin{algorithmic}[1]
\State $MaxDisc = 1$
\State $SolutionSets = \{\}$
\State $ApproxSets = \{\}$ 
\State $PriorSets = \{\}$
\State $MixRatioSets = \{\}$
\State $TauSets = \{\}$
\State $DiscCntr = 0$
\State estimate $\sigma_{1}^{2}$, ..., $\sigma_{N_{l}}^{2}$ from the sample means approximating the standard deviations of $e_{1}$, ..., $e_{N_{l}}$. \vspace{0.05cm}  
\State $\delta^{2} := \mathrm{max}\big\{\sigma_{1}^{2}, ..., \sigma_{N_{l}}^{2}\big\}$
\State $\boldsymbol{\Sigma} := \delta^{-2} \cdot \mathrm{diag}\big(\sigma_{1}^{2}, ..., \sigma_{N_{l}}^{2}\big)$
\State $p_{i} = \frac{i - 1}{N_{frac} - 1}, \; i = 1, ..., N_{frac}$ \label{MixingRatios}
\State $N_{frac} = 201$
\State $N_{mean} = 5$
\State $I_{min} = \{\}$
     \For{$i = 1 \; \textbf{to} \; N_{frac}$} \label{BeginKernelMatrices}
	\For{$j = 1 \; \textbf{to} \; N_{l}$}
	\State compute $k_{tot}(l_{j}) + i n_{tot}(l_{j})$ from $k_{1}(l_{j}) + i n_{1}(l_{j})$ and $k_{2}(l_{j}) + i n_{2}(l_{j})$ 
 using \eqref{Lorentz-Lorenz} \hspace*{1.0cm} with $f_{1} = p_{i}$ and $f_{2} = 1 - p_{i}$
	\EndFor
	\For{$k = 1 \; \textbf{to} \; m$}
	\State compute kernel matrix $\boldsymbol{K}_{ik}$ for $p_{i}$ and the collocation grid with 
           $N_{k}$ points \hspace*{1.0cm} using $k_{tot}(l_{1}) + i n_{tot}(l_{1}), ..., k_{tot}(l_{N_{l}}) + i n_{tot}(l_{N_{l}})$
	\EndFor
     \EndFor \label{EndKernelMatrices}
     \For{$k = 1 \; \textbf{to} \; m$} \label{BeginMainLoop}
           \State $S_{k} = \{\}$
           \State $A_{k} = \{\}$
           \State $P_{k} = \{\}$
           \State $M_{k} = \{\}$
           \State $T_{k} = \{\}$
	\State $R = \{\}$
	\State $RM_{min} = \infty$
%
%
	\For{$i = 1 \; \textbf{to} \; N_{frac}$} \label{BeginScan}
	\State $\boldsymbol{n}_{lsqnng} = \underset{\boldsymbol{n} \: \in \: \mathbb{R}^{N_{k}}}{\mathrm{argmin}} \; \textstyle{\frac{1}{2}}\|\boldsymbol{\Sigma}^{-\frac{1}{2}}(\boldsymbol{K}_{ik}\boldsymbol{n} - \boldsymbol{e}_{real})\|_2^{2} \;\; \text{s.t.} \; \boldsymbol{n} \geq 0$ \vspace{0.1cm}
	\State $R = R \cup \big\{\|\boldsymbol{\Sigma}^{-\frac{1}{2}}(\boldsymbol{K}_{ik}\boldsymbol{n}_{lsqnng} - \boldsymbol{e}_{real})\|_2^{2}\big\}$
          \EndFor
	\For{$i = 1 \; \textbf{to} \; N_{frac} - N_{mean} + 1$}
	\State $RM =\mathrm{mean}(R(i), R(i + 1), ..., R(i + N_{mean} - 1))$
                \If{$RM < RM_{min}$}
	     \State $RM_{min} = RM$	
	     \State $t_{min} = \{i, i + 1, ..., i + N_{mean} - 1\}$
	     \EndIf
	\EndFor \label{EndScan}
	\State $t_{cur} = \{t_{min}(1), t_{min}(3), t_{min}(5)\}$ \label{model_inds}
	\State $N_{cur} = |t_{cur}|$
	\For{$i = 1 \; \textbf{to} \; N_{cur}$} \label{BeginModelGeneration}
                 \For{$j=0 \; \textbf{to} \; s$}
      \algstore{myalg3}
      \end{algorithmic}
   \end{algorithm}  

\newpage

   \begin{algorithm}
      \begin{algorithmic}
         \algrestore{myalg3} 
		\If{$R(t_{cur}(i)) < \tau_{j}N_{l}\delta^{2} \; \land \; \tau_{j}N_{l}\delta^{2} < \|\boldsymbol{\Sigma}^{-\frac{1}{2}}\boldsymbol{e}_{real}\|_2^{2}$} \vspace{0.1cm} \label{DiscPrincipleCheck2}
	      	\State $\text{compute }\gamma_{kij}\text{ such that}$ \vspace{0.1cm}	
	      	\State $\boldsymbol{n}_{trial} = \underset{\boldsymbol{n} \: \in \: \mathbb{R}^{N_{i}}}{\mathrm{argmin}} \; \textstyle{\frac{1}{2}}\|\boldsymbol{\Sigma}^{-\frac{1}{2}}(\boldsymbol{K}_{t_{cur}(i), k}\boldsymbol{n} - \boldsymbol{e}_{real})\|_2^{2} + \textstyle{\frac{1}{2}}\gamma_{kij}\boldsymbol{n}^{T}\boldsymbol{R}_{k}\boldsymbol{n} \;\; \text{s.t.} \; \boldsymbol{n} \geq 0$ 
 	     	\State $\text{with } \|\boldsymbol{\Sigma}^{-\frac{1}{2}}(\boldsymbol{K}_{t_{cur}(i), k}\boldsymbol{n}_{trial} - \boldsymbol{e}_{real})\|_2^{2} = \tau_{j}N_{l}\delta^{2}$	\vspace{0.05cm}  
		\EndIf
	  	\If{$\boldsymbol{n}_{trial}\text{ exists}$} \label{BeginContainers}
		\State $S_{k} = S_{k} \cup \{\boldsymbol{n}_{trial}\}$
		\State $A_{k} = A_{k} \cup \{\boldsymbol{K}_{t_{cur}(i), k}\}$
		\State $P_{k} = P_{k} \cup \{\gamma_{kij}\boldsymbol{R}_{k}\}$
		\State $M_{k} = M_{k} \cup \{p_{t_{cur}(i)}\}$
		\State $T_{k} = T_{k} \cup \{\tau_{j}\}$
		\EndIf \label{EndContainers}
                 \EndFor
	  \EndFor \label{EndModelGeneration}
	  \If{$S_{k}, \; A_{k}, \; P_{k}, \; M_{k} \text{ and } T_{k} \text{ not empty}$} 
	  \State $SolutionSets = SolutionSets \cup \{S_{k}\}$
	  \State $ApproxSets = ApproxSets \cup \{A_{k}\}$
	  \State $PriorSets = PriorSets \cup \{P_{k}\}$
             \State $MixRatioSets = MixRatioSets \cup \{M_{k}\}$
	  \State $TauSets = TauSets \cup \{T_{k}\}$
	  \State $DiscCntr = DiscCntr + 1$
	  \EndIf  
	  \If{$DiscCntr == MaxDisc$}
	  \State $\text{break}$ \label{ExitLoop}
	  \EndIf
     \EndFor
\end{algorithmic}
\end{algorithm}

In line \ref{MixingRatios} the aerosol fraction parameter interval $[0, 1]$ is approximated with a linearly spaced grid. For each discrete aerosol fraction $p_{i}$ the approximation $\boldsymbol{K}_{ik}$ to the linear operator $K_{p_{i}}$ is computed in lines \ref{BeginKernelMatrices} to \ref{EndKernelMatrices} for all model space orders $N_{k}$.

In line \ref{BeginMainLoop} the main loop for the model generation begins. Note that we first run through all model orders from $1$ to $m$ beginning with the coarsest models before we iterate through all aerosol fractions $p_{i}$. This means that we perform the residual-based search strategy motivated in Proposition \ref{unique_fraction} for each model space separately, where we start with the coarsest model and refine it if necessary.

In lines \ref{BeginScan}-\ref{EndScan} the residuals of the unregularized reconstructions are calculated, and a scan to find the minimal mean of $N_{mean}$ solutions corresponding to successive parameters $p_{i}, p_{i+1}, ..., p_{i + N_{mean} - 1}$ is performed. A subset of the indices $i, i+1, ..., i + N_{mean} - 1$ corresponding to the residuals with minimal mean is selected in line \ref{model_inds}. By filtering out some of the models corresponding to the parameters $p_{i}, p_{i+1}, ..., p_{i + N_{mean} - 1}$ with small residuals we ensure that the models to be compared are not too similar. The selected indices are used for the actual model generation in lines \ref{BeginModelGeneration}-\ref{EndModelGeneration}. Here we loop through all preselected Morozov safety parameters $\tau_{1}, ..., \tau_{s}$ and we propose with them the possible residual values $\tau_{j}N_{l}$ for the discrepancy principle. In line \ref{DiscPrincipleCheck2} it is checked if the discrepancy principle is applicable. 

If the model generation step is successful, the obtained reconstructions accompanied by their kernel and regularization matrices and their aerosol fraction and residual parameters are stored in the containers $S_{k}$, $A_{k}$, $P_{k}$, $M_{k}$ and $T_{k}$ in lines \ref{BeginContainers}-\ref{EndContainers}.

Finally if the model generation is successful for $MaxDisc$ model spaces, the model generation loop is terminated in line \ref{ExitLoop}.

\subsection{Model Selection for Two-Component Aerosols under Nonnegativity Constraints}

Not only the model generation procedure has to be generalized to the case of a two-component aerosol, but also the model-selection framework presented in Section \ref{section_model_selection} needs to be generalized as well. Here we are not just comparing models with different model spaces but also with different underlying operators $K_{p}$. Thus prior probabilities are also needed for the parameters $p$ which determine the linear operators---or more precisely their approximations---to be compared. Let $k$ label the model  dimensions $N_{k}$, Let $i$ run through the indices for the aerosol-fraction parameters $p_{i}$, where $i$ depends on $k$, and let $j$ run through all Morozov safety parameters $\tau_{j}$ used for the model generation, where $j$ depends on $k$ and $i$. Then we can compute the model posterior probabilities by
\begin{equation}
p(N_{k}, \boldsymbol{K}_{ik}, \gamma_{kij} | \boldsymbol{e}) =  \frac{p(\boldsymbol{e} | N_{k}, \boldsymbol{K}_{ik}, \gamma_{kij})p( N_{k}, \boldsymbol{K}_{ik}, \gamma_{kij})}{\sum_{u}\sum_{v(u)}\sum_{w(u,v)}p(\boldsymbol{e} | N_{u}, \boldsymbol{K}_{vu}, \gamma_{uvw})p( N_{u}, \boldsymbol{K}_{vu}, \gamma_{uvw)})}
\end{equation}
We assume that $p(\boldsymbol{K}_{vu})$ and $p(N_{u}, \gamma_{uvw})$ are independent and thus
\begin{equation*}
p( N_{u}, \boldsymbol{K}_{vu}, \gamma_{uvw}) = p(N_{u}, \gamma_{uvw})p(\boldsymbol{K}_{vu}).
\end{equation*}
We select $p(\boldsymbol{K}_{vu})$ to be uniform and adopt $p(N_{u}, \gamma_{uvw})$ from Section \ref{subsetion_model_selection}. This leads to
\begin{equation}
p( N_{u}, \boldsymbol{K}_{vu}, \gamma_{uvw}) = \frac{1}{N_{total}},
\end{equation}
where $N_{total}$ is the total number of triplets $\big(u, \; v(u), \; w(u,v)\big)$.

Then the model-selection algorithm proceeds in the same way as Algorithm \ref{model_selection}, so we do not restate here. The differences to Section \ref{subsetion_model_selection} are that we have already set $MaxDisc = 1$ in the model generation step and that the single container $A_{1}$ stores kernel matrices approximating different operators $K_{p_{i}}$. While in principle the algorithm could continue to compare different discretizations, this only lead to worse results in our simulations. Therefore, once the algorithm finds a discretization level for which reconstructions are at all possible for any of the safety factors, we stop the refinement and simply focus on the problem of identifying the volume fraction.

\section{Numerical Results for Two-Component Aerosols}

\subsection{Numerical Study}

We conducted a numerical study of our inversion algorithm with almost the same settings as the last section but extended for the retrieval of volume fractions of a two-component aerosol. We used the same wavelength grid as in Sections \ref{simulated_measurements} and simulated the same model size distributions as in Section \ref{model_size_distributions}. We selected air as ambient medium as well. We extended the grid of Morozov safety parameters to
\begin{equation*}
\tau_{1} = 0.5, \tau_{2} = 0.6, ..., \tau_{16} = 2.0.
\end{equation*}
If when running through all model spaces none of these safety factors yielded a solution, we performed in this extreme case an another run of the model generation step using a second grid of safety factors given by
\begin{equation*}
\tau_{1} = 2.5, \tau_{2} = 3.0, ..., \tau_{6} = 5.0.
\end{equation*}

This time we did not just simulate an original aerosol consisting purely of $\mathrm{H}_{2}\mathrm{O}$ but instead generated with \eqref{Lorentz-Lorenz} refractive indices of $\mathrm{H}_{2}\mathrm{O}$ and $\mathrm{CsI}$ mixtures for the scattering particles. Here the volume fractions of $\mathrm{H}_{2}\mathrm{O}$ ranged through a set of preselected percentages, namely 
\begin{equation*}
0, \; 11, \; 22, \; 33, \; 44, \; 56, \; 67, \; 78, \; 89 \; \text{ and } \; 100.
\end{equation*}
For each of the $100$ parameters for the log-normal, RRSB or Hedrih distributions we also now have the above $10$ fractions.  This results in a total of $1000$ cases to simulate.

As preparation to run Algorithm \ref{model_generation_2comp} we computed the kernel matrices depending on the water volume fraction parameter $p \in [0, 1]$ for $p_{i} = \frac{i - 1}{100}$, $i = 1, ..., 101$ and interpolated each kernel matrix entry with a cubic spline on a linearly spaced grid with $201$ points covering $[0, 1]$ to increase further the resolution in $p$. Thus we have $N_{frac} = 201$ in Algorithm \ref{model_generation_2comp}.

Another important difference to Section \ref{simulated_measurements} is that the noise level was taken to be only $5$\% of the true measurement values $e(l_{i})$ instead of $30$\%. Thus the noisy measurement data vector $\boldsymbol{e}$ was here modeled with
\begin{equation*}
(\boldsymbol{e})_{i} = e(l_{i}) + \delta_{i} \quad \text{ with } \; \delta_{i} \sim \mathcal{N}(0,(0.05 \cdot e(l_{i}))^{2}), \quad i = 1, ..., N_{l}.
\end{equation*}
We had to take this lower value because the problem of retrieving the aerosol fractions additionally to the size distributions is much more ill-posed than simply reconstructing the size distribution when the scattering material is known.

To investigate the quality of the reconstructions we computed their $L^{2}$-errors relative to the original size distribution. We list them separately for each of the ten original water fractions. We proceed this way for all of our simulation results.

Furthermore we determined the deviations of the reconstructed water volume fractions from the original ones, e.g.\ when the original fraction was $22$\% and $p_{recon} \in [0, 1]$ the retrieved fraction parameter, we calculated the deviation by $| 22 - 100 \cdot p_{recon} | \%$. This showed us how well one can investigate the unknown two-component aerosol only from FASP measurements using our extended inversion algorithm.

We also report how often the inversions failed. There were two main reasons for inversion failures: the first when the relative $L^{2}$ was greater than or equal to $100$\%, the second when the fraction deviation was greater than or equal to $50$\%. In both cases the reconstruction cannot give any reasonable information about the true size distribution and the true scattering material anymore. Note that in our simulations we returned by default $\boldsymbol{n} \equiv 0$ and $p_{recon} = 0.5$ when no reconstruction could be found in any of the model spaces. For brevity we only list those original fractions where inversion failures occurred

Finally we list the average and worst case inversion run times over all $1000$ simulations. 

\subsubsection{Average $L^{2}$-Errors}

\begin{table}[h!]
\centering
\begin{tabular}{||m{2.55cm}||m{2.5cm}|m{2.5cm}|m{2.5cm}||}
\hline
\multicolumn{4}{||c||}{\multirow{3}{*}{}} \tabularnewline[-2.5ex]
\multicolumn{4}{||c||}{Log-Normal Distribution} \tabularnewline
\multicolumn{4}{||c||}{} \tabularnewline[-2.5ex]
\hline
\multirow{4}{*}{\vspace{0.9cm} original water} & \multicolumn{3}{c||}{} \tabularnewline[-2.5ex]
& \multicolumn{3}{c||}{average $L^{2}$-errors (\%)} \tabularnewline
& \multicolumn{3}{c||}{} \tabularnewline[-2.5ex]
\cline{2-4}
& \centering \multirow{2}{*}{Tikhonov} & \centering min. first & \centering \multirow{2}{*}{Twomey} \tabularnewline
\vspace{-0.6cm} volume percent & & \centering fin. diff. & \tabularnewline
\hline
\hline
0 \% & \centering 33.5018 & \centering 33.3939 & \centering 33.6631 \tabularnewline
\hline
11 \% & \centering 30.0691 & \centering 30.1394 & \centering 31.1483 \tabularnewline
\hline
22 \% & \centering 28.8329 & \centering 28.5679 & \centering 29.6568 \tabularnewline
\hline
33 \% & \centering 24.9249 & \centering 24.9295 & \centering 25.8386 \tabularnewline
\hline
44 \% & \centering 23.4809 & \centering 23.4363 & \centering 24.6004 \tabularnewline
\hline
56 \% & \centering 22.6134 & \centering 22.2024 & \centering 23.2278 \tabularnewline
\hline
67 \% & \centering 20.4546 & \centering 20.2646 & \centering 20.3727 \tabularnewline
\hline
78 \% & \centering 18.9590 & \centering 18.7913 & \centering 19.3874 \tabularnewline
\hline
89 \% & \centering 18.5494 & \centering 18.7126 & \centering 18.9227 \tabularnewline
\hline
100 \% & \centering 17.9452 & \centering 17.7468 & \centering 18.6161 \tabularnewline
\hline
\end{tabular}
\end{table}

\newpage

\begin{table}[h!]
\centering
\begin{tabular}{||m{2.55cm}||m{2.5cm}|m{2.5cm}|m{2.5cm}||}
\hline
\multicolumn{4}{||c||}{\multirow{3}{*}{}} \tabularnewline[-2.5ex]
\multicolumn{4}{||c||}{RRSB Distribution} \tabularnewline
\multicolumn{4}{||c||}{} \tabularnewline[-2.5ex]
\hline
\multirow{4}{*}{\vspace{0.9cm} original water} & \multicolumn{3}{c||}{} \tabularnewline[-2.5ex]
& \multicolumn{3}{c||}{average $L^{2}$-errors (\%)} \tabularnewline
& \multicolumn{3}{c||}{} \tabularnewline[-2.5ex]
\cline{2-4}
& \centering \multirow{2}{*}{Tikhonov} & \centering min. first & \centering \multirow{2}{*}{Twomey} \tabularnewline
\vspace{-0.6cm} volume percent & & \centering fin. diff. & \tabularnewline
\hline
\hline
0 \% & \centering 38.4939 & \centering 37.8096 & \centering 36.4143 \tabularnewline
\hline
11 \% & \centering 31.7809 & \centering 31.2368 & \centering 30.5024 \tabularnewline
\hline
22 \% & \centering 30.4269 & \centering 29.9565 & \centering 29.5171 \tabularnewline
\hline
33 \% & \centering 27.7925 & \centering 27.4651 & \centering 26.2830 \tabularnewline
\hline
44 \% & \centering 24.2137 & \centering 24.1607 & \centering 23.3831 \tabularnewline
\hline
56 \% & \centering 21.5396 & \centering 21.8053 & \centering 21.3910 \tabularnewline
\hline
67 \% & \centering 19.5181 & \centering 19.8905 & \centering 19.5506 \tabularnewline
\hline
78 \% & \centering 16.9927 & \centering 16.7956 & \centering 16.7368 \tabularnewline
\hline
89 \% & \centering 14.4005 & \centering 14.4332 & \centering 14.2989 \tabularnewline
\hline
100 \% & \centering 12.0441 & \centering 11.5948 & \centering 11.6963 \tabularnewline
\hline
\end{tabular}
\end{table}

\begin{table}[h!]
\centering
\begin{tabular}{||m{2.55cm}||m{2.5cm}|m{2.5cm}|m{2.5cm}||}
\hline
\multicolumn{4}{||c||}{\multirow{3}{*}{}} \tabularnewline[-2.5ex]
\multicolumn{4}{||c||}{Hedrih Distribution} \tabularnewline
\multicolumn{4}{||c||}{} \tabularnewline[-2.5ex]
\hline
\multirow{4}{*}{\vspace{0.9cm} original water} & \multicolumn{3}{c||}{} \tabularnewline[-2.5ex]
& \multicolumn{3}{c||}{average $L^{2}$-errors (\%)} \tabularnewline
& \multicolumn{3}{c||}{} \tabularnewline[-2.5ex]
\cline{2-4}
& \centering \multirow{2}{*}{Tikhonov} & \centering min. first & \centering \multirow{2}{*}{Twomey} \tabularnewline
\vspace{-0.6cm} volume percent & & \centering fin. diff. & \tabularnewline
\hline
\hline
0 \% & \centering 17.2492 & \centering 17.0990 & \centering 16.5398 \tabularnewline
\hline
11 \% & \centering 17.5805 & \centering 17.5246 & \centering 17.1164 \tabularnewline
\hline
22 \% & \centering 17.3665 & \centering 17.1217 & \centering 16.5678 \tabularnewline
\hline
33 \% & \centering 15.0919 & \centering 15.0406 & \centering 14.6545 \tabularnewline
\hline
44 \% & \centering 16.9073 & \centering 16.8601 & \centering 16.8571 \tabularnewline
\hline
56 \% & \centering 15.5153 & \centering 15.5509 & \centering 15.2833 \tabularnewline
\hline
67 \% & \centering 14.5155 & \centering 14.4545 & \centering 14.5169 \tabularnewline
\hline
78 \% & \centering 16.7670 & \centering 16.7708 & \centering 16.6434 \tabularnewline
\hline
89 \% & \centering 13.3015 & \centering 13.3173 & \centering 13.2861 \tabularnewline
\hline
100 \% & \centering 11.3083 & \centering 11.2952 & \centering 11.4170 \tabularnewline
\hline
\end{tabular}
\end{table}

\newpage

\subsubsection{Average Water Fraction Deviation}

\begin{table}[h!]
\centering
\begin{tabular}{||m{2.55cm}||m{2.5cm}|m{2.5cm}|m{2.5cm}||}
\hline
\multicolumn{4}{||c||}{\multirow{3}{*}{}} \tabularnewline[-2.5ex]
\multicolumn{4}{||c||}{Log-Normal Distribution} \tabularnewline
\multicolumn{4}{||c||}{} \tabularnewline[-2.5ex]
\hline
\multirow{4}{*}{\vspace{0.9cm} original water} & \multicolumn{3}{c||}{} \tabularnewline[-2.5ex]
& \multicolumn{3}{c||}{average water fraction deviation (\%)} \tabularnewline
& \multicolumn{3}{c||}{} \tabularnewline[-2.5ex]
\cline{2-4}
& \centering \multirow{2}{*}{Tikhonov} & \centering min. first & \centering \multirow{2}{*}{Twomey} \tabularnewline
\vspace{-0.6cm} volume percent & & \centering fin. diff. & \tabularnewline
\hline
\hline
0 \% & \centering 10.7150 & \centering 10.7750 & \centering 10.9350 \tabularnewline
\hline
11 \% & \centering 7.3700 & \centering 7.4400 & \centering 7.5000 \tabularnewline
\hline
22 \% & \centering 6.1750 & \centering 6.0450 & \centering 6.1050 \tabularnewline
\hline
33 \% & \centering 4.3700 & \centering 4.3400 & \centering 4.4800 \tabularnewline
\hline
44 \% & \centering 3.9550 & \centering 3.9750 & \centering 4.0350 \tabularnewline
\hline
56 \% & \centering 3.2000 & \centering 3.1900 & \centering 3.3600 \tabularnewline
\hline
67 \% & \centering 2.6050 & \centering 2.5650 & \centering 2.5850 \tabularnewline
\hline
78 \% & \centering 2.2350 & \centering 2.2150 & \centering 2.1950 \tabularnewline
\hline
89 \% & \centering 2.2000 & \centering 2.2300 & \centering 2.2100 \tabularnewline
\hline
100 \% & \centering 1.4150 & \centering 1.2950 & \centering 1.4150 \tabularnewline
\hline
\end{tabular}
\end{table}

\begin{table}[h!]
\centering
\begin{tabular}{||m{2.55cm}||m{2.5cm}|m{2.5cm}|m{2.5cm}||}
\hline
\multicolumn{4}{||c||}{\multirow{3}{*}{}} \tabularnewline[-2.5ex]
\multicolumn{4}{||c||}{RRSB Distribution} \tabularnewline
\multicolumn{4}{||c||}{} \tabularnewline[-2.5ex]
\hline
\multirow{4}{*}{\vspace{0.9cm} original water} & \multicolumn{3}{c||}{} \tabularnewline[-2.5ex]
& \multicolumn{3}{c||}{average water fraction deviation (\%)} \tabularnewline
& \multicolumn{3}{c||}{} \tabularnewline[-2.5ex]
\cline{2-4}
& \centering \multirow{2}{*}{Tikhonov} & \centering min. first & \centering \multirow{2}{*}{Twomey} \tabularnewline
\vspace{-0.6cm} volume percent & & \centering fin. diff. & \tabularnewline
\hline
\hline
0 \% & \centering 6.5800 & \centering 6.6800 & \centering 6.6600 \tabularnewline
\hline
11 \% & \centering 5.3000 & \centering 5.2900 & \centering 5.2400 \tabularnewline
\hline
22 \% & \centering 4.7100 & \centering 4.6700 & \centering 4.6300 \tabularnewline
\hline
33 \% & \centering 3.7300 & \centering 3.7800 & \centering 3.8000 \tabularnewline
\hline
44 \% & \centering 3.7300 & \centering 3.6900 & \centering 3.6500 \tabularnewline
\hline
56 \% & \centering 3.1650 & \centering 3.1950 & \centering 3.1450 \tabularnewline
\hline
67 \% & \centering 2.3650 & \centering 2.3850 & \centering 2.3450 \tabularnewline
\hline
78 \% & \centering 1.8050 & \centering 1.8550 & \centering 1.8950 \tabularnewline
\hline
89 \% & \centering 1.3750 & \centering 1.4150 & \centering 1.3850 \tabularnewline
\hline
100 \% & \centering 0.4650 & \centering 0.4050 & \centering 0.4050 \tabularnewline
\hline
\end{tabular}
\end{table}

\newpage

\begin{table}[h!]
\centering
\begin{tabular}{||m{2.55cm}||m{2.5cm}|m{2.5cm}|m{2.5cm}||}
\hline
\multicolumn{4}{||c||}{\multirow{3}{*}{}} \tabularnewline[-2.5ex]
\multicolumn{4}{||c||}{Hedrih Distribution} \tabularnewline
\multicolumn{4}{||c||}{} \tabularnewline[-2.5ex]
\hline
\multirow{4}{*}{\vspace{0.9cm} original water} & \multicolumn{3}{c||}{} \tabularnewline[-2.5ex]
& \multicolumn{3}{c||}{average water fraction deviation (\%)} \tabularnewline
& \multicolumn{3}{c||}{} \tabularnewline[-2.5ex]
\cline{2-4}
& \centering \multirow{2}{*}{Tikhonov} & \centering min. first & \centering \multirow{2}{*}{Twomey} \tabularnewline
\vspace{-0.6cm} volume percent & & \centering fin. diff. & \tabularnewline
\hline
\hline
0 \% & \centering 5.0800 & \centering 5.0000 & \centering 5.0600 \tabularnewline
\hline
11 \% & \centering 6.3050 & \centering 6.3150 & \centering 6.1850 \tabularnewline
\hline
22 \% & \centering 4.2600 & \centering 4.4500 & \centering 4.1000 \tabularnewline
\hline
33 \% & \centering 4.4350 & \centering 4.4550 & \centering 4.5350 \tabularnewline
\hline
44 \% & \centering 3.2100 & \centering 3.2300 & \centering 3.0800 \tabularnewline
\hline
56 \% & \centering 3.3750 & \centering 3.3950 & \centering 3.2450 \tabularnewline
\hline
67 \% & \centering 1.8700 & \centering 1.9000 & \centering 2.0000 \tabularnewline
\hline
78 \% & \centering 2.7750 & \centering 2.7550 & \centering 2.7850 \tabularnewline
\hline
89 \% & \centering 1.9150 & \centering 1.9150 & \centering 1.8750 \tabularnewline
\hline
100 \% & \centering 1.0050 & \centering 0.9850 & \centering 1.0250 \tabularnewline
\hline
\end{tabular}
\end{table}

\subsection{Average Model Space Dimensions}

\begin{table}[h!]
\centering
\begin{tabular}{|m{2.5cm}|m{2.5cm}|m{2.5cm}|}
\cline{1-3}
\multicolumn{3}{|c|}{\multirow{3}{*}{}} \tabularnewline[-2.5ex]
\multicolumn{3}{|c|}{Log-Normal Distribution} \tabularnewline
\multicolumn{3}{|c|}{} \tabularnewline[-2.5ex]
\cline{1-3}
\multicolumn{3}{|c|}{} \tabularnewline[-2.5ex]
\multicolumn{3}{|c|}{average model space dimensions} \tabularnewline
\multicolumn{3}{|c|}{} \tabularnewline[-2.5ex]
\cline{1-3}
\centering \multirow{2}{*}{Tikhonov} & \centering min. first & \centering \multirow{2}{*}{Twomey} \tabularnewline
& \centering fin. diff. & \tabularnewline
\hline
\centering 7.5550 & \centering 7.5550 & \centering 7.5550 \tabularnewline
\hline
\end{tabular}
\end{table}

\begin{table}[h!]
\centering
\begin{tabular}{|m{2.5cm}|m{2.5cm}|m{2.5cm}|}
\cline{1-3}
\multicolumn{3}{|c|}{\multirow{3}{*}{}} \tabularnewline[-2.5ex]
\multicolumn{3}{|c|}{RRSB Distribution} \tabularnewline
\multicolumn{3}{|c|}{} \tabularnewline[-2.5ex]
\cline{1-3}
\multicolumn{3}{|c|}{} \tabularnewline[-2.5ex]
\multicolumn{3}{|c|}{average model space dimensions} \tabularnewline
\multicolumn{3}{|c|}{} \tabularnewline[-2.5ex]
\cline{1-3}
\centering \multirow{2}{*}{Tikhonov} & \centering min. first & \centering \multirow{2}{*}{Twomey} \tabularnewline
& \centering fin. diff. & \tabularnewline
\hline
\centering 11.3940 & \centering 11.3940 & \centering 11.3940 \tabularnewline
\hline
\end{tabular}
\end{table}

\begin{table}[h!]
\centering
\begin{tabular}{|m{2.5cm}|m{2.5cm}|m{2.5cm}|}
\cline{1-3}
\multicolumn{3}{|c|}{\multirow{3}{*}{}} \tabularnewline[-2.5ex]
\multicolumn{3}{|c|}{Hedrih Distribution} \tabularnewline
\multicolumn{3}{|c|}{} \tabularnewline[-2.5ex]
\cline{1-3}
\multicolumn{3}{|c|}{} \tabularnewline[-2.5ex]
\multicolumn{3}{|c|}{average model space dimensions} \tabularnewline
\multicolumn{3}{|c|}{} \tabularnewline[-2.5ex]
\cline{1-3}
\centering \multirow{2}{*}{Tikhonov} & \centering min. first & \centering \multirow{2}{*}{Twomey} \tabularnewline
& \centering fin. diff. & \tabularnewline
\hline
\centering 6.3330 & \centering 6.3330 & \centering 6.3330 \tabularnewline
\hline
\end{tabular}
\end{table}

\subsection{Extreme Cases}

When the deviation of the retrieved aerosol fraction from the true one exceeded $50\%$ or the $L^{2}$-error between reconstruction and true solution was bigger than $100\%$ we had to regard the reconstruction as failed. We now list when these failures occurred.
There were no failed simulations with the Hedrih distribution.

\subsubsection{Reconstruction Failures}

\begin{table}[h!]
\centering
\begin{tabular}{||m{2.55cm}||m{2.5cm}|m{2.5cm}|m{2.5cm}||}
\hline
\multicolumn{4}{||c||}{\multirow{3}{*}{}} \tabularnewline[-2.5ex]
\multicolumn{4}{||c||}{Log-Normal Distribution} \tabularnewline
\multicolumn{4}{||c||}{} \tabularnewline[-2.5ex]
\hline
\multirow{4}{*}{\vspace{0.9cm} original water} & \multicolumn{3}{c||}{} \tabularnewline[-2.5ex]
& \multicolumn{3}{c||}{number of $L^{2}$-errors $\geq 100$ \% (out of 100)} \tabularnewline
& \multicolumn{3}{c||}{} \tabularnewline[-2.5ex]
\cline{2-4}
& \centering \multirow{2}{*}{Tikhonov} & \centering min. first & \centering \multirow{2}{*}{Twomey} \tabularnewline
\vspace{-0.6cm} volume percent & & \centering fin. diff. & \tabularnewline
\hline
\hline
0 \% & \centering 1 & \centering 1 & \centering 1 \tabularnewline
\hline
11 \% & \centering 1 & \centering 1 & \centering 1 \tabularnewline
\hline
44 \% & \centering 0 & \centering 0 & \centering 1 \tabularnewline
\hline
\end{tabular}
\end{table}

\begin{table}[h!]
\centering
\begin{tabular}{||m{2.55cm}||m{2.5cm}|m{2.5cm}|m{2.5cm}||}
\hline
\multicolumn{4}{||c||}{\multirow{3}{*}{}} \tabularnewline[-2.5ex]
\multicolumn{4}{||c||}{RRSB Distribution} \tabularnewline
\multicolumn{4}{||c||}{} \tabularnewline[-2.5ex]
\hline
\multirow{4}{*}{\vspace{0.9cm} original water} & \multicolumn{3}{c||}{} \tabularnewline[-2.5ex]
& \multicolumn{3}{c||}{number of $L^{2}$-errors $\geq 100$ \% (out of 100)} \tabularnewline
& \multicolumn{3}{c||}{} \tabularnewline[-2.5ex]
\cline{2-4}
& \centering \multirow{2}{*}{Tikhonov} & \centering min. first & \centering \multirow{2}{*}{Twomey} \tabularnewline
\vspace{-0.6cm} volume percent & & \centering fin. diff. & \tabularnewline
\hline
\hline
0 \% & \centering 8 & \centering 8 & \centering 7 \tabularnewline
\hline
11 \% & \centering 1 & \centering 1 & \centering 1 \tabularnewline
\hline
22 \% & \centering 1 & \centering 1 & \centering 1 \tabularnewline
\hline
33 \% & \centering 1 & \centering 1 & \centering 1 \tabularnewline
\hline
44 \% & \centering 1 & \centering 0 & \centering 0 \tabularnewline
\hline
67 \% & \centering 1 & \centering 1 & \centering 1 \tabularnewline
\hline
\end{tabular}
\end{table}

\subsubsection{Water-Fraction Retrieval Failure}

\begin{table}[h!]
\centering
\begin{tabular}{||m{2.55cm}||m{2.5cm}|m{2.5cm}|m{2.5cm}||}
\hline
\multicolumn{4}{||c||}{\multirow{3}{*}{}} \tabularnewline[-2.5ex]
\multicolumn{4}{||c||}{Log-Normal Distribution} \tabularnewline
\multicolumn{4}{||c||}{} \tabularnewline[-2.5ex]
\hline
\multirow{4}{*}{\vspace{0.9cm} original water} & \multicolumn{3}{c||}{} \tabularnewline[-2.5ex]
& \multicolumn{3}{c||}{number of deviations $\geq 50$ \% (out of 100)} \tabularnewline
& \multicolumn{3}{c||}{} \tabularnewline[-2.5ex]
\cline{2-4}
& \centering \multirow{2}{*}{Tikhonov} & \centering min. first & \centering \multirow{2}{*}{Twomey} \tabularnewline
\vspace{-0.6cm} volume percent & & \centering fin. diff. & \tabularnewline
\hline
\hline
0 \% & \centering 3 & \centering 3 & \centering 3 \tabularnewline
\hline
11 \% & \centering 3 & \centering 3 & \centering 3 \tabularnewline
\hline
22 \% & \centering 1 & \centering 1 & \centering 1 \tabularnewline
\hline
\end{tabular}
\end{table}

\begin{table}[h!]
\centering
\begin{tabular}{||m{2.55cm}||m{2.5cm}|m{2.5cm}|m{2.5cm}||}
\hline
\multicolumn{4}{||c||}{\multirow{3}{*}{}} \tabularnewline[-2.5ex]
\multicolumn{4}{||c||}{RRSB Distribution} \tabularnewline
\multicolumn{4}{||c||}{} \tabularnewline[-2.5ex]
\hline
\multirow{4}{*}{\vspace{0.9cm} original water} & \multicolumn{3}{c||}{} \tabularnewline[-2.5ex]
& \multicolumn{3}{c||}{number of deviations $\geq 50$ \% (out of 100)} \tabularnewline
& \multicolumn{3}{c||}{} \tabularnewline[-2.5ex]
\cline{2-4}
& \centering \multirow{2}{*}{Tikhonov} & \centering min. first & \centering \multirow{2}{*}{Twomey} \tabularnewline
\vspace{-0.6cm} volume percent & & \centering fin. diff. & \tabularnewline
\hline
\hline
0 \% & \centering 1 & \centering 1 & \centering 1 \tabularnewline
\hline
\end{tabular}
\end{table}

\newpage

\subsubsection{Average and Worst-Case Run Times}

\begin{table}[h!]
\centering
\begin{tabular}{m{2.0cm}|m{2.5cm}|m{2.5cm}|m{2.5cm}|}
\cline{2-4}
& \multicolumn{3}{c|}{\multirow{3}{*}{}} \tabularnewline[-2.5ex]
& \multicolumn{3}{c|}{Log-Normal Distribution} \tabularnewline
& \multicolumn{3}{c|}{} \tabularnewline[-2.5ex]
\cline{2-4}
& \multicolumn{3}{c|}{} \tabularnewline[-2.5ex]
& \multicolumn{3}{c|}{run times (s)} \tabularnewline
& \multicolumn{3}{c|}{} \tabularnewline[-2.5ex]
\cline{2-4}
& \centering \multirow{2}{*}{Tikhonov} & \centering min. first & \centering \multirow{2}{*}{Twomey} \tabularnewline
& & \centering fin. diff. & \tabularnewline
\hline
\multicolumn{1}{|c|}{average} & \centering 1.9700 & \centering 2.0189 & \centering 2.0453 \tabularnewline
\hline
\multicolumn{1}{|c|}{worst case} & \centering 8.5281 & \centering 7.7737 & \centering 7.8988 \tabularnewline
\hline
\end{tabular}
\end{table}

\begin{table}[h!]
\centering
\begin{tabular}{m{2.0cm}|m{2.5cm}|m{2.5cm}|m{2.5cm}|}
\cline{2-4}
& \multicolumn{3}{c|}{\multirow{3}{*}{}} \tabularnewline[-2.5ex]
& \multicolumn{3}{c|}{RRSB Distribution} \tabularnewline
& \multicolumn{3}{c|}{} \tabularnewline[-2.5ex]
\cline{2-4}
& \multicolumn{3}{c|}{} \tabularnewline[-2.5ex]
& \multicolumn{3}{c|}{run times (s)} \tabularnewline
& \multicolumn{3}{c|}{} \tabularnewline[-2.5ex]
\cline{2-4}
& \centering \multirow{2}{*}{Tikhonov} & \centering min. first & \centering \multirow{2}{*}{Twomey} \tabularnewline
& & \centering fin. diff. & \tabularnewline
\hline
\multicolumn{1}{|c|}{average} & \centering 2.4997 & \centering 2.5488 & \centering 2.5536 \tabularnewline
\hline
\multicolumn{1}{|c|}{worst case} & \centering 16.6366 & \centering 16.6889 & \centering 16.5122 \tabularnewline
\hline
\end{tabular}
\end{table}

\begin{table}[h!]
\centering
\begin{tabular}{m{2.0cm}|m{2.5cm}|m{2.5cm}|m{2.5cm}|}
\cline{2-4}
& \multicolumn{3}{c|}{\multirow{3}{*}{}} \tabularnewline[-2.5ex]
& \multicolumn{3}{c|}{Hedrih Distribution} \tabularnewline
& \multicolumn{3}{c|}{} \tabularnewline[-2.5ex]
\cline{2-4}
& \multicolumn{3}{c|}{} \tabularnewline[-2.5ex]
& \multicolumn{3}{c|}{run times (s)} \tabularnewline
& \multicolumn{3}{c|}{} \tabularnewline[-2.5ex]
\cline{2-4}
& \centering \multirow{2}{*}{Tikhonov} & \centering min. first & \centering \multirow{2}{*}{Twomey} \tabularnewline
& & \centering fin. diff. & \tabularnewline
\hline
\multicolumn{1}{|c|}{average} & \centering 1.4207 & \centering 1.4429 & \centering 1.4629 \tabularnewline
\hline
\multicolumn{1}{|c|}{worst case} & \centering 4.7495 & \centering 4.9839 & \centering 5.0766 \tabularnewline
\hline
\end{tabular}
\end{table}

\subsection{Conclusion}

A common trend in the results is that the average $L^{2}$-error decreases with increasing original water volume fraction. For log-normal distributions it averages around $31$\% for water fractions ranging from $0$ to $22$\% and for RRSB distributions it decreases from around $37$ to $30$\% for the same fractions. This poor behavior can also be seen in the numbers of reconstruction failures, which only occurred for log-normal and RRSB distributions and mostly for water fractions below or equal to $44$\%. For higher fractions ranging from $56$ to $100$\% the average $L^{2}$ was always below $22$\% for all three size distribution classes and even improving towards $100$\%.

The water volume fractions deviations behaved in a similar way. They decreased for increasing original fractions, which means that the quality of the water fraction retrieval was improving towards higher original fractions. Water fraction retrieval failures only happened for log-normal and RRSB distributions, when the original fractions were below or equal to $22$\%. Again the differences in the deviations depending on the priors were only marginal.

The worst case run times never succeeded our thirty-seconds limit. Even in the extreme cases they stayed below $17$ seconds. The average run times ranged from ca. $1.7$ to $3.3$ seconds.

We can conclude that with the settings made in previous section the analysis of two-component aerosols is possible satisfying our demands on run time and accuracy. The standard deviations of the noise in all single measurements have to be reduced from $30\%$ of the true extinction values to $5\%$ in order to obtain results of comparable quality as in Section \ref{num_study}. 

\section{Outlook}

The model selection problem, i.e.\ to select appropriate model spaces $\mathbb{R}^{N}$, can also be treated with \textit{Markov Chain Monte Carlo Methods}. Here the posterior distribution is defined as multidimensional distribution living on all model spaces and sampled by the Monte Carlo method. We plan to compare our methods developed to these.
 
\section{Acknowledgement}

We thank Prof. Dr. Hans-Josef Allelein and Bj\"{o}rn Krupa for providing us with experimental FASP measurement data, which was the guide for our simulated measurement data.

\begin{figure}[htbp]
\begin{minipage}{0.75\textwidth}
This work is sponsored by the German Federal Ministry of Education and Research (BMBF) under the contract number 02NUK022.

Responsibility for the content of this report lies with the authors.
\end{minipage}
\hfill
\begin{minipage}{0.225\textwidth}
\includegraphics[width=\textwidth]{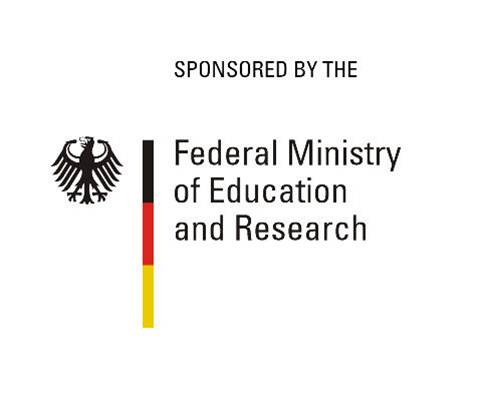}
\end{minipage}
\end{figure}

\bibliographystyle{ieeetr}
\bibliography{literature}

\end{document}